\documentclass[10pt]{article}
\usepackage{amsmath}
\usepackage{amssymb,amsmath,mathrsfs}
\usepackage{color}
\usepackage{graphicx}
\usepackage{pstricks}
\usepackage{epsfig,subfigure}
\let\pa=\partial

\let\s=\sigma
\let\f=\frac
\let\p=\psi
\let\om=\omega
\let\G= \Gamma
\let\D=\Delta

\let\Om=\Omega
\let\wt=\widetilde

\def\cA{{\cal A}}

\def\cD{{\cal D}}

\def\cH{{\cal H}}

\def\cN{{\cal N}}

\def\cV{{\cal V}}

\def\no{\noindent}
\def\na{\nabla}

\def\p{\partial}

\def\dv{\mbox{div}}

\def\curl{\mathop{\rm curl}\nolimits}

\def\C{\mathop{\bf C\kern 0pt}\nolimits}
\def\DD{\mathop{\bf D\kern 0pt}\nolimits}
\def\K{\mathop{\bf K\kern 0pt}\nolimits}
\def\N{\mathop{\bf N\kern 0pt}\nolimits}
\def\Q{\mathop{\bf Q\kern 0pt}\nolimits}
\def\R{\mathop{\bf R\kern 0pt}\nolimits}
\def\T{\mathop{\bf T\kern 0pt}\nolimits}

\newcommand{\ef}{ \hfill $ \blacksquare $ \vskip 3mm}

\newcommand{\beq}{\begin{equation}}
\newcommand{\eeq}{\end{equation}}
\newcommand{\ben}{\begin{eqnarray}}
\newcommand{\een}{\end{eqnarray}}
\newcommand{\beno}{\begin{eqnarray*}}
\newcommand{\eeno}{\end{eqnarray*}}

\newtheorem{rmk}{Remark}[section]

\newtheorem{prop}{Proposition}[section]
\renewcommand{\theequation}{\thesection.\arabic{equation}}

\newtheorem{theorem}{Theorem}[section]

\newtheorem{lemma}[theorem]{Lemma}
\newtheorem{proposition}[theorem]{Proposition}

\begin{document}
\title{ Water waves problem with surface tension in a corner domain I: A priori estimates with constrained contact angle}
\author{Mei Ming$^1$ and Chao Wang$^2$ \\[2mm]
{\small $ ^1$ School of Mathematics, Sun Yat-sen University, Guangzhou 510275, China
  }\\[2mm]
{\small E-mail: mingm@mail.sysu.edu.cn}\\[2mm]
{\small $ ^2$ School of  Mathematical Science, Peking University, Beijing 100871, China}\\[2mm]
{\small E-mail: wangchao@math.pku.edu.cn}
}
\date{ }
\maketitle

\begin{abstract}
We study the two-dimensional water waves problem with surface tension in the case when there is a non-zero contact angle between the free surface  and the bottom. In the presence of surface tension, dissipations take place at the contact point. Moreover, when the contact angle is less than $\pi /6$, no singularity appears in our settings. Using elliptic estimates in corner domains and a geometric approach, we prove  an a priori estimate for the water waves problem. 
\end{abstract}

\renewcommand{\theequation}{\thesection.\arabic{equation}}
\setcounter{equation}{0}

\tableofcontents

\section{Introduction}
In this paper, we consider the  water waves problem in an  unbounded two-dimensional corner domain $\Om_t$ with an  upper free surface $\Gamma_t$ and a fixed flat bottom $\Gamma_b$.  Specifically, if we parametrize our domain  by denoting $\G_t=\{(x,z)\,|\, z=\eta(t, x)\}$ and $\G_b=\{(x,z)\,|\, z=b(x)\}$,   we can  write at time $t$ that
\beno
\Om_t= \{(x,z)\,|\,b(x)<z<\eta(t,x) \}.
\eeno
To avoid  technical complexity,  we  assume that $\G_t$ and $\G_b$ only have one intersection point  $X_c$ (the contact point) at the left end (See figure 1). The contact angle between the free surface and the bottom is denoted by $\om(t)$ (or sometimes simply $\om$).  

Without loss of generality, the contact point  is at the origin when $t=0$. Moreover, we assume that the bottom $\G_b$ is a line segment near the origin, and $\G_b$ becomes a horizontal line away from the origin.
Our domain $\Om_t$ has a  finite depth,  which means that there exists a constant $h>0$ such that the distance between $\G_t,\,\G_b$ is always less that $h$. 
One can tell that $\Om_t$ corresponds to the scene of sea waves moving near a beach,  and the contact point denotes the intersection point between the sea and the rigid bottom in two dimensional case.

\medskip

The water waves problem investigates an ideal fluid with a free surface, which is supposed  to be inviscid and incompressible.  We assume that  the fluid is under the influence of  gravity in the domain $\Om_t$ and  surface tension on the free surface $\G_t$. Moreover, the problem is also assumed to be irrotational.

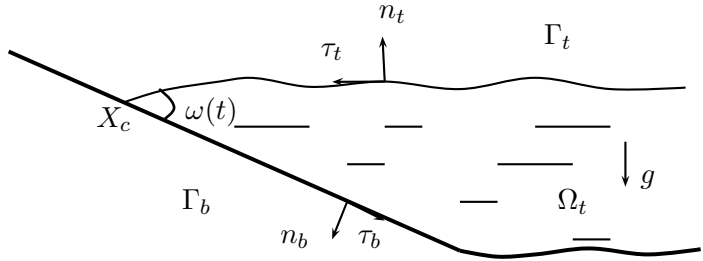
\begin{figure}
\begin{pspicture}(-5,-3)(5,2.5)
\pscurve[showpoints=false](-0.45,-0.17)(0,-0.03)(0.8,0.1)(1.2,0.15)(2,0.03)(3,0.1)(4,0)(5,0.15)(6,0)(7,0.02)
\psline[linewidth=1.5pt](-2,0.5)(4,-2.15)
\pscurve[linewidth=1.5pt](4,-2.15)(4.5,-2.23)(5,-2.2)(5.7,-2.12)(6.5,-2.2)(7.2,-2.13)
\rput(3.1,1){$ n_t$}
\psline{->}(3,0.1)(2.97,0.7)
\rput(2.3,0.5){$\tau_t$}
\psline{->}(3,0.1)(2.3,0.096)
\rput(1.8,-2){$ n_b$}
\psline{->}(2.5,-1.5)(2.3,-2)
\rput(2.8,-2){$\tau_b$}
\psline{->}(2.5,-1.5)(3,-1.75)
\rput(6.5,-1.2){$ g$}
\psline{->}(6.2,-.7)(6.2,-1.3)
\psline(1,-0.5)(2,-0.5)
\psline(3,-0.5)(3.5,-0.5)
\psline(5,-0.5)(6,-0.5)
\psline(2.5,-1)(3,-1)
\psline(4.5,-1)(5.5,-1)
\psline(4,-1.5)(4.5,-1.5)
\psline(5.5,-2)(6,-2)
\pscurve[linewidth=1pt](0,0)(0.2, -0.25)(0.07,-0.4)
\rput(-0.6,-0.4){$X_c$}
\rput(0.7, -0.3){$\om(t)$}
\rput(5.5,-1.5){$\Om_t$}
\rput(5.3,0.7){$\G_t$}
\rput(0.5,-1.5){$\G_b$}
\end{pspicture}
\caption{The corner domain}
\label{figdomain}
\end{figure}


The water waves problem in our case involves surface tension on the free surface. Compared to the case without surface tension, it is important to find a proper condition at the contact point.   
In fact, T. Young \cite{Young} had proved a long time ago (1805, Young's law) that, in the stationary case, the (stationary) contact angle $\om_s$ is a physical constant depending only on the materials of the bottom and the fluid:
\[
\cos{\om_s}=\f{[\gamma]}{\sigma}
\]
where $\sigma$ is the surface tension coefficient along the free surface and $[\gamma]=\gamma_{1}-\gamma_{2}$ with $\gamma_{1},\gamma_{2}$  measures of the free-energy per unit length associated to the solid-vapor and solid-fluid interaction respectively.

Based on this  theory, when the fluid is moving, 
a modified condition from W. Ren and W. E \cite{RE} is adopted in our paper, which takes the fluid speed at the contact point into account:
\beq\label{condition at X_c}
\beta_c v_c=\sigma(\cos{\om_s}-\cos{\om}) \qquad\hbox{or equivalently}\quad \beta_cv_c=[\gamma]-\sigma \cos\om,
\eeq
where 
\[
v_c=-v\cdot \tau_b
\] 
denotes the fluid speed at the contact point along the bottom and $\beta_c$ is the effective friction coefficient determined by interfacial widths, interactions between the fluid and the bottom, and the normal stress contributions.

In fact, this condition tells us that slip velocity is dominated by the unbalanced Young stress,   which  is an effective variation of Young's law  for  stationary contact angles \cite{Young}.
This kind of conditions are very common and widely discussed, see \cite{BEIMR, CDA, SA, GL}.

\medskip

Based on the condition above  and letting $v$ be  the fluid velocity and $P$ be the total pressure, we consider the following water waves problem  in  the domain $\Om_t$ at time $t>0$:
 \beq\label{WW}
\mbox{(WW)}\qquad \left\{
\begin{array}{l}
\pa_t v+v\cdot \na v=-\na P-g, \\
\dv v=0,\quad \curl v=0,\qquad \hbox{on}\quad \Om_t\\
P|_{\Gamma_t}=\sigma \kappa,\\
\pa_t+v\cdot \na \quad \textrm{is tangent to}\quad \Gamma_t,\\
v\cdot n_b |_{\Gamma_b}=0,\\
\beta_c v_c=\sigma(\cos{\om_s}-\cos{\om}),\qquad\hbox{at}\quad X_c
\end{array}
\right.
\eeq
with  $\kappa$  the mean curvature of the free surface and $g$ the constant gravity vector.

\medskip
Before stating our result, firstly we  would like to recall  some previous works on the well-posedness of classical water waves problems. Here `classical' means water waves in smooth domains, where smooth domains refer to domains with smooth boundaries.  Compared to  smooth domains,  when we say non-smooth domains, we always refer to domains with Lipschitz boundaries. For example, our domain $\Om_t$ with a corner is a non-smooth domain.

In the case when surface tension is ignored, some early works such as V.I. Nalimov \cite{Na}, H. Yosihara \cite{Yo1}, and W. Craig \cite{Craig} established   the local well-posedness with small data  in 2 dimensional case. The  local well-posedness of general initial data for 2 and 3 dimensional was solved by S. Wu \cite{Wu1, Wu2} in the case of infinite depth.  D. Lannes \cite{Lannes} considered the water waves problem in the case of finite depth under Eulerian settings.  Moreover, H. Lindblad \cite{Lin} proved the existence of solutions for the general problem of a liquid body in vacuum and P. Zhang and Z. Zhang \cite{ZZ} showed the  local well-posedness for the rotational problem.

Concerning the problem with surface tension, H. Yoisihara \cite{Yo2} proved an early result on the local well-posedness with small data in infinite-depth case. T. Iguchi \cite{Iguchi} and D. Ambrose \cite{Ambrose} studied the local well-posedness of the irrotational problem in 2 dimensional case. B. Schweizer \cite{Sch} showed the existence for the general 3 dimensional problem. 
In \cite{SZ, SZ2}, J. Shatah and C. Zeng  proved  a priori estimates and the local well-posedness (even when the fluid domains are not simply connected),  where they used a geometric approach  and the fluid was rotational and with surface tension. A similar geometric approach had also been used by K. Beyer and M. G\"unther\cite{BG1, BG2} to study the irrotational problem for star-shaped domains.  
Moreover, D. Coutand and S. Shkoller \cite{CS} proved local well-posedness for the rotational problem under Lagrangian coordinates.  

For the global well-posedness for the water waves, the first result was given by S. Wu \cite{Wu4} who proved  the almost global existence for the gravity problem in two dimensions. Later, P. Germain, N. Masmoudi and J. Shatah \cite{GMS} and S. Wu \cite{Wu5} proved the global existence of gravity waves in three dimensions respectively.  Moreover, T. Alazard and  J.M. Delort \cite{AD} and A.D. Ionescu and F. Pusateri \cite{IP}  studied  global regularity for the gravitational water waves system in two dimensions independently. Recently, J. Hunter, M. Ifrim and D. Tataru \cite{HIT1,HIT2, HIT3} used the conform mapping method to give another proof of the global existence for the gravitational problem in two dimensions. 
  
There are more works on water waves, and we only mention some of them here: Ambrose and Masmoudi \cite{AM1, AM2},  M. Ogawa and A. Tani \cite{OT1, OT2}, B. Alvarez-Samaniego and D. Lannes \cite{AL},  T. Alazard, Burq and Zuily \cite{ABZ}, M. Ming and Z. Zhang \cite{MZ}, M. Ming, P. Zhang and Z. Zhang \cite{MZZ} e.t.c..

\medskip
Secondly, we  recall some recent results on water waves  in non-smooth domains.  In fact, theoretical research on the non-smooth domain only started several years ago and there remains a lot of open problems.  There was a work by T. Alazard, N. Burq and C.  Zuily \cite{ABZ15} for some right angle with vertical walls when there is no surface tension,  where they used  symmetrizing and periodizing  to turn this problem into a classical water waves problem.  Later,  R.H. Kinsey and S. Wu \cite{WuK} and  S. Wu \cite{Wu3} proved a priori estimates and the local well-posedness  for the two dimensional water waves with angled crests, where a conformal mapping is used to convert the boundary singularities. T. de Poyferr\'e \cite{Poyferre} gave  a priori estimates for the rotational water waves problem in a compact domain (corresponds to a beach type) in general $n$ dimensions, where the contact angle is smaller than a dimensional constant and therefore no singularity appears. Moreover, this work is done in absence of  surface tension. On the other hand, a very recent work \cite{LM} by D. Lannes and G. M\'etivier solved the local well-posedness for the Green-Naghdi equations in a beach-type domain, which is an approximation model for the water waves under shallow-water regime.

\medskip
In our paper, we consider a priori estimates of the water waves problem $\mbox{(WW)}$ in a beach-type domain $\Om_t$, which is two dimensional with only one contact point. Different from T. de Poyferr\'e \cite{Poyferre}, we take  surface tension on the free surface  into account, which leads to a big difference in the whole energy formulation and will be discussed very soon.  

To study the water waves problem $\mbox{(WW)}$, the first main difficulty comes from the Dirichlet-Neumann operator or equivalently the related elliptic systems on corner domains. As already explained in details in our previous work \cite{MW},  non-smooth domain generates singularities from related elliptic systems.  
Moreover, Remark 5.20 \cite{MW} tells us directly that smaller contact angle leads to higher regularities for elliptic systems. 
Motived from this point, we want to  study the water waves problem under a proper formulation firstly {\it with no singularity}. In fact,  to avoid singularity, we  will work under a constrained contact angle $\om$ and a comparatively low regularity (will be explained in the following main theorem) as the first step for our project.   Meanwhile, we adopt  the geometric approach from  J. Shatah and C. Zeng \cite{SZ, SZ2}, which turns out to be a good choice for corner domains.
  
   Compared to  the classical water waves, our problem is also variational, but the variation formulation is very different. In fact, one can see in Section 3 that, the Lagrangian Action contains a potential at the contact point, and meanwhile there is also a dissipation related to the contact point in the variation equation \eqref{dissipation eqn}. In a word,  our variation formulation is new compared to the classical water waves.

Now we stress the role of  surface tension. As mentioned above, the energy formulations with and without surface tension are completely different. When there is no surface tension, the friction at the contact point is ignored as in \cite{Poyferre}, so the energy is conserved. On the other hand, when surface tension is taken into considerations,  condition \eqref{condition at X_c} appears, and system $\mbox{(WW)}$ generates a different basic energy in Eulerian coordinates
\beq\label{basic energy}
E_0=\f12\displaystyle\int_{\Om_t}|v|^2dX+g\displaystyle\int_{\Om_t}X\cdot e_z dX+\sigma S(u)+[\gamma]X\cdot \tau_b|_{X_c},
\eeq  
where the last term denotes the friction potential at the contact point.
Moreover, one can find in Section 4 that system $\mbox{(WW)}$ satisfies the following  dissipation equation:
\beq\label{dissipation eqn} 
\f d{dt} E_0+\beta_c|v_c|^2\big|_{X_c}=0
\eeq 
with $\beta_c>0$ the friction coefficient and $v_c$ the speed at the contact point along the bottom. 
So one can see that  some energy dissipation takes place at the contact point.  This dissipation formulation is indeed similar as that in Y. Guo and I. Tice \cite{GT}, which works on the contact line problem for Stokes equation.

As a result, the dissipation leads to a big difference in the a priori estimate: The terms involving the contact point need to be treated carefully, which is a completely new part in water waves problem. Similarly as in \cite{SZ}, our energy estimate is proved firstly for the equation of the main part $J$ from $\na P$, and then we go back to $\mbox{(WW)}$. Compared to \cite{SZ}, another difference in the estimates relates to the bottom $\G_b$, which also needs much care since our estimates are performed very often in variational sense. Consequently,  some special Sobolev spaces from P. Grisvard \cite{PG1} such as $\tilde H^{\f12}(\G_b),\,\tilde H^{-\f12}(\G_b)$ are used in our paper.

\medskip
Now it's the time to state the main theorem. To begin with, we introduce the energy functional
\[
E(t)=\|\na_{\tau_t}J^\bot\|^2_{L^2(\Gamma_t)} +\|\cD_t J\|^2_{L^2(\Om_t)}+\|\Gamma_t\|^2_{H^{\f52}}+\|v\|^2_{L^2(\Om_t)},
\] 
and the dissipation at the contact point
\[
F(t)=\big|(\sin \om)\na_{{\bf \tau}_t}J^\perp |_{X_c}\big|^2.
\]
The main theorem is presented here.
\begin{theorem}\label{main theorem}
Assume that the initial data $(\G_0,v_0) \in H^4\times H^3(\Om_0)$ and the initial contact angle $\om_0\in(0,\f\pi6)$.
Let $(\G_t, v) \in H^4\times C\big(H^3(\Om_t)\big)$ be a  solution of $\mbox{(WW)}$, then there exists a constant $T_0$ depending on the initial data such that the following a priori estimate holds
\beno
\sup_{0\leq t\leq T_0}E(t)+\int_0^{T_0} F(t)dt\leq E(0)+ \int_0^{T_0} P(E(t))dt, 
\eeno
where $P(\cdot)$ is a polynomial with positive constant coefficients depending on $\sigma,\beta_c, [\gamma], \G_b$.
\end{theorem}

\begin{rmk}
{\it We consider the irrotational case in this paper, but our formulation may also work for the rotational case. 
}
\end{rmk}

\begin{rmk}
{\it In this paper, we need at most $H^4$ estimates for related elliptic systems in $\Om_t$. One can see directly from Remark 5.20 \cite{MW} that, to avoid the singularity, one needs the contact angle $\om\in(0,\f \pi 6)$.  The energy estimates for the water waves problem with a  general angle remains an open problem.
}
\end{rmk}

\begin{rmk}
{\it 
Notice that the dissipation $F(t)$ contains $\sin \om$, which means smaller contact angle $\om$ leads to smaller dissipation as long as $\na_{\tau_t}J^\perp|_{X_c}$ remains the same. On the other hand, we have $J\in H^\f32(\Om_t)$ in our settings, so the dissipation term $\int_0^{T_0} F(t)dt$ on the left side of the a priori estimate is some kind of smoothing estimates such that $\na_{\tau_t}J^\perp$ makes sense at the contact point $X_c$ for $t\in [0,T_0]$. 
}
\end{rmk}
 

\no {\bf Organization of this paper.} Section 2 introduces some notations used in this paper. In Section 3, we explain the problem in a geometric approach and prove that it is variational. In Section 4, a dissipation equation is deduced. In Section 5, we  recall some trace theorems and elliptic estimates from our previous paper. Section 6 deals with some commutators and derive the equation for $J$. In the end, we prove the a priori estimates in Section 7.

\bigskip

\section{Notations}
\setcounter{equation}{0}

- $\Om_0$ is the initial domain at time $t=0$, and $\Om_t$ is the domain at time $t$.\\
- We denote by $Y$ a point in $\Om_0$, and by $X$ a point in $\Om_t$.\\
- $X_c$ is the coordinate of the contact point at time $t$, which corresponds to $Y_c\in \Om_0$ satisfying $X_c=u(t,Y_c)$.\\ 
- $v_c=-v\cdot \tau_b|_{X_c}$ is the speed of the contact point $X_c$ along the bottom $\G_b$.\\
- We denote by $\tilde n_b,\,\tilde \tau_b$ unit orthogonal extensions onto $\Om_t$ for $n_b,\,\tau_b$ on $\G_b$.\\
- $D_t=\p_t+\na_v$ is the material derivative.\\
- $M^*$ denotes the transport of a matrix $M$.\\
- $A\cdot B$ denotes the inner product of two vectors or two matrices $A,B$.\\
- $w^\perp$ on $\G_t$:   $w\cdot n_t$ for a vector $w\in T_X\G_t$.  \\
- $w^\top$ on $\G_t$: $(w\cdot \tau_t)\,\tau_t$. Sometimes we also use $w^\top$ on $\G_b$ with a similar definition.\\
- $\Pi$: the second fundamental form where $\Pi(w)=\na_w n_t\in T_X\G_t$ for $w\in T_X\G_t$.\\
- $\Pi(v,w)$ denotes $\Pi(v)\cdot w$. Moreover, $\Pi$ is symmetric: $\Pi(v,w)=\Pi(w,v)$.\\
- $|\Pi|^2=tr (\Pi \,\Pi^*)$.\\
- $\kappa=tr \Pi=\na_{\tau_t}n_t\cdot \tau_t$ is the mean curvature.\\
- $\bar \cD$ is the covariant derivative on $\G$, and $\cD$ is $\bar \cD$ in Eulerian coordinates.\\
- $(\cD\cdot \Pi)(w)=(\cD_{\tau_t}\Pi)(w)\cdot \tau_t=\cD_{\tau_t}(\Pi(w))-\Pi(\cD_{\tau_t}w)$.\\
- $\D_{\G_t}$ is the Beltrami-Laplace operator on $\G_t$:
\[
\D_{\G_t}f=\cD^2f(\tau_t,\tau_t)=\cD\cdot (\na^\top f)=\na_{\tau_t}\na_{\tau_t}f-\na_{\cD_{\tau_t}\tau_t}f.
\]
- $\cD^2f(\tau_1, \tau_2)=D^2f(\tau_1, \tau_2)-(\Pi(\tau_1)\cdot\tau_2)\na_{n_t} f$  for any two vectors $\tau_1,\tau_2$.\\
- $\cH(f)$ or $f_\cH$ is the harmonic extension for some function $f$ on $\G_t$, which is defined by the elliptic system
\[
\left\{\begin{array}{ll}
\D \cH(f)=0,\qquad\hbox{on}\quad \Om_t,\\
\cH(f)|_{\G_t}=f,\quad \na_{n_b}\cH(f)|_{\G_b}=0.
\end{array}\right.
\]
- $\D^{-1}(h,g)$ denotes  the solution of the system
\[
\left\{\begin{array}{ll}
\D u=h\qquad \hbox{on}\quad \Om_t\\
u|_{\G_t}=0,\qquad \na_{n_b}u|_{\G_b}=g.
\end{array}
\right.
\]
- $[\gamma]=\gamma_1-\gamma_2$, where $\gamma_1,\,\gamma_2$ are the surface tension coefficients denoting the solid-air and solid-fluid interactions respectively.  \\
- $\beta_c$ is the effective friction coefficient determined by interfacial widths, interactions between the fluid and the bottom, and the normal stress contributions.\\
- $g$ denotes the constant gravity vector or the gravity coefficient.\\
- $P(E(t))$: Some polynomial for the energy $E(t)$ with positive constant coefficients.\\
- The Sobolev norm $H^s$ on $\Om_t$ can be defined by  restrictions
\beno
\|u\|_{H^{s}(\Om)}=\inf\{\|U\|_{H^s(\R^2)}, \, U|_{\Om}=u\}.
\eeno
- $\tilde H^{\f12}(\G_b)$ is a subspace of $H^{\f12}(\G_b)$ related to corner domains
\[
\tilde H^{\f12}(\G_b)=\Big\{u\in \dot{H}^{\f12}(\G_b)\Big| \,\rho^{-\f12}u\in L^2(\G_b)\Big\}\] where $\dot H^\f12(\G_b)$ is the closure of $\mathscr D(\G_b)$ in $H^s(\G_b)$,  and $\rho=\rho(X)$ is the distance (arc length) between the point $X\in \G_b$ and the left end $X_c$.  The norm is defined as 
\[
\|u\|^2_{\tilde H^{\f12}}=\|u\|^2_{H^\f12}+\int_{\G_b} \rho^{-1}|u|^2dX.
\] 
Moreover,  we use $\tilde H^{-\f12}(\G_b)$ to denote the dual space of $\tilde H^\f12(\G_b)$. For more details, see \cite{PG1}. \\

\bigskip

\section{Geometry and variation}
\setcounter{equation}{0}

In this section, we introduce the geometry behind the water waves problem following the notations from \cite{SZ}. One can see that our problem is also variational in nature, while the new point here is about the dissipation at the contact point. In the end, we will compute the second variation of the basic energy $E_0$ to find out the leading-order term in the linearization of system $\mbox{(WW)}$, which turns out to be the same as that in \cite{SZ}.

Let $X=u(t,Y)$ for any  $Y\in \Omega_0$ be the Lagrangian  coordinates map solving 
\[
\f{dX}{dt}=v(t,X),\quad X(0)=Y.
\] 
So the velocity can be denoted as $v=u_t\circ u^{-1}$. 

Since $v$ is divergence free, the trajectory map $u$ is volume-preserving.
Accordingly, we can define the manifold
\[
\Gamma=\{\Phi:\Om_0\rightarrow \R^2|\quad\Phi \ \hbox{is a volume-perserving homeomorphism}\},
\] 
and consequently the tangent space of $\G$ is given by
\[
T_\Phi \G=\{\bar w: \Om_0\rightarrow \R^2|\, w=\bar w\circ \Phi^{-1}\ \hbox{satisfying}\ \na\cdot w=0\, \ \hbox{on} \ \Phi(\Om_0)\ \hbox{and}\  w\cdot n_b|_{\G_b}=0 \}
\] 
where $\G_b$, $n_b$ denote the bottom and the unit outward normal vector of $\Phi(\Om_0)$ respectively. 

Based on the tangent space, we  also need to consider $(T_\Phi \G)^\perp$ and the Hodge decomposition. In fact, for any vector field $w: \Phi(\Om_0)\rightarrow \R^2$, we have the Hodge decomposition 
\[
w=w_1-\na q
\] 
where $\bar w_1=w_1\circ \Phi\in T_\Phi \G$ and $q\in (T_\Phi \G)^\perp$  decided by $\Phi,w$:
\[
\left\{\begin{array}{ll}
-\D q=\na\cdot w\qquad \hbox{on}\quad \Om_t\\
q|_{\G_t}=0,\qquad \na_{n_b}q|_{\G_b}=-w\cdot n_b.
\end{array}
\right.
\]
Moreover, one can see that 
\[
(T_\Phi \G)^\perp=\{-(\na q)\circ \Phi\  \big|\   q|_{\G_t}=0\}
\]
with $\G_t$ the upper surface of $\Phi(\Om_0)$.

Now we are ready to apply the Hodge decomposition.  For  a path $u(t,\cdot)\in \G$ with $\bar v=u_t$, and any vector $\bar w(t,\cdot)\in T_{u(t)}\G$, we decompose $\bar w_t$ to have the covariant derivative $\bar \cD_t\bar w$ and the second fundamental form $II_{u(t)}(\bar w,\bar v)$ satisfying
\[
\bar w_t=\bar\cD_t\bar w+II_{u(t)}(\bar w,\bar v),
\] 
where 
\[
\bar\cD_t\bar w\in T_{u(t)}\G,\quad II_{u(t)}(\bar w,\bar v)\in (T_{u(t)}\G)^\perp.
\] From the Hodge decomposition we know that
\[
II_{u(t)}(\bar w,\bar v)=-(\na P_{w,v})\circ u
\] 
with $P_{w,v}$ solving
\[
\left\{\begin{array}{ll}
-\D P_{w,v}=tr (\na w\na v),\qquad \hbox{on}\quad \Om_t\\
P_{w,v}|_{\G_t}=0,\quad \na_{n_b}P_{w,v}|_{\G_b}=w\cdot \na_vn_b|_{\G_b}.
\end{array}
\right.
\] 
As a result, denoting $\cD_t w=(\bar \cD_t\bar w)\circ u^{-1}$, we rewrite the decomposition as
\[
\cD_t w=D_t w+\na P_{w,v}.
\]

\medskip
Before proceeding further, we recall the following computation directly from \cite{SZ}:
\[
D_tn_t=-\big((\na v)^*n_t\big)^\top\qquad\hbox{on}\quad \G_t.
\] 
Moreover, since our domain is two dimensional in this paper, we know directly that 
$\tau_t$ is the parallel-transporting tangent basis satisfies
\[
D_t\tau_t=(\na_{\tau_t}v\cdot n_t)n_t\qquad\hbox{and}\quad \cD_{\tau_t}\tau_t=0\qquad\hbox{on}\quad \G_t
\]  where $\cD$ is Eulerian-coordinates version of the covariant derivative $\bar \cD$ on $\G$. 
These expressions will be used repeatedly in the following sections.

\subsection {Lagrangian variation}
Our problem has a variation formulation with dissipations at the corner. In fact, we introduce the Lagrangian Action on a time interval $[0,T]$:
\[
I(u)=\displaystyle\int^T_0\Big(\displaystyle\int_{\Om_t}\f12|v|^2dX-g\displaystyle\int_{\Om_t}zdX-\sigma S(u)-[\gamma] X\cdot \tau_b|_{X_c}\Big)
\] where $\tau_b$ near $X_c$ is a constant vector due to the definition of the bottom, and the surface potential
\beq\label{S(u) def}
S(u)=\displaystyle\lim_{A\rightarrow +\infty}\Big(\displaystyle\int_{\G_t\cap\{x\le A\}}ds-\displaystyle\int_{\G_*\cap\{x\le A\}}ds\Big)
\eeq
with $\G_*$ some reference surface.

On the other hand, we  define the dissipation at the contact point
\[
F(u,q)=\f12\beta_c\big|(q\circ u^{-1})\cdot \tau_b|_{X_c}\big|^2
\] 
with $q=u_t$.

For any path $u(s,t,\cdot)\in \G$ with $\bar w=u_s|_{s=0}$ and $w=\bar w\circ u^{-1}|_{s=0}$ satisfying 
\[
\bar w|_{t=0}=\bar w|_{t=T}=0,
\] 
we will show that the Euler equation and the condition \eqref{condition at X_c} from system $\mbox{(WW)}$  can be deduced from the variation formulation 
\[
\langle I'(u), \bar w\rangle=\displaystyle\int^T_0\langle F_q(u,u_t),\bar w\rangle dt.
\] 
To prove this, we firstly start with the left side and compute on $\Om_0$ to find:
\[
\langle I'(u), \bar w\rangle=\displaystyle\int^T_0\Big(\displaystyle\int_{\Om_0}u_t\cdot \bar w_tdY-g\displaystyle\int_{\Om_0}e_z\cdot \bar wdY-\sigma \langle S'(u),\bar w\rangle-[\gamma]\bar w\cdot \tau_b|_{Y_c}\Big)dt.
\] 
A direct calculation as in \cite{SZ} leads to
\[
\langle S'(u),\bar w \rangle=\displaystyle\int_{\G_t}(\kappa w^\bot+\cD\cdot w^\top)ds,
\] 
where in our 2 dimensional case, we can write in particular that
\[
\cD\cdot w^\top=\na_{\tau_t}(w\cdot \tau_t)=-\f d{ds}(w\cdot \tau_t)
\] with $s$ the arclength parameter on $\G_t$ starting from $X_c$. Consequently, we find 
\[
\displaystyle\int_{\G_t}\cD\cdot w^\top ds=-\displaystyle\int^\infty_{0}\f d{ds}(w\cdot \tau_t)ds=w\cdot \tau_t|_{X_c},
\] 
which infers that
\[
\langle S'(u),\bar w\rangle=\displaystyle\int_{\G_t}\kappa w^\perp ds+w\cdot \tau_t|_{X_c}
=\displaystyle\int_{\Om_t}\na\kappa_\cH\cdot wdX+w\cdot \tau_t|_{X_c}.
\]
Notice that compared to the classical case, there is an extra term concerning the contact point. 
As a result, we arrive at
\[
\langle I'(u), \bar w\rangle=\displaystyle\int^T_0\Big(-\displaystyle\int_{\Om_t}(D_t v+g e_z+\sigma\na \kappa_{\cH})\cdot w dX-(\sigma w\cdot \tau_t+[\gamma]w\cdot \tau_b)|_{X_c}\Big)dt.
\]
Secondly, the right side of the variation formulation turns out to be
\[
\displaystyle\int^T_0\langle F_q(u,u_t),\bar w\rangle dt =\displaystyle\int^T_0\beta_c(v\cdot \tau_b)\,(w\cdot \tau_b)|_{X_c}dt,
\]
which together with the left side implies the following equality: 
\[
\displaystyle\int^T_0\Big(-\displaystyle\int_{\Om_t}(D_t v+g e_z+\sigma\na \kappa_{\cH})\cdot w dX-(\sigma w\cdot \tau_t+[\gamma]w\cdot \tau_b)|_{X_c}\Big)dt=
\displaystyle\int^T_0\beta_c(v\cdot \tau_b)\,(w\cdot \tau_b)|_{X_c}dt.
\]
Consequently, we retrieve the Euler equation 
\[
D_t v=-\na P_{v,v}-\sigma\na \kappa_{\cH}-g e_z\qquad\hbox{on}\quad \Om_t
\] 
where the total pressure $P$ in $\mbox{(WW)}$ is decomposed into two parts
\[
P=P_{v,v}+\sigma\kappa_\cH
\]
with $P_{v,v}$  defined by the Hodge decomposition and $\sigma\kappa_\cH$ from the mean curvature $\kappa$.

Moreover, we also retrieve the equality at the corner:
\[
-\displaystyle\int^T_0(\sigma w\cdot \tau_t+[\gamma]w\cdot \tau_b)|_{X_c}dt=
\displaystyle\int^T_0\beta_c(v\cdot \tau_b)\,(w\cdot \tau_b)|_{X_c}dt.
\] Remembering the notations
\[
v\cdot \tau_b=-v_c,\quad w\cdot \tau_b=-w_c,
\]  
one can have
\[
w\cdot \tau_t=w_c\cos\om.
\] 
Substituting these computations into the  equality above, we  derive the condition at the contact point:
\[
\beta_c v_c=[\gamma]-\sigma \cos\om,
\] which can also be written as 
\[
\beta_cv_c=\sigma(\cos{\om_s}-\cos \om).
\] Here recall that $\om_s$ stands for  the  static contact angle $\om_s$.

\medskip

\subsection{Second variation of the basic energy}
Recall that the basic energy $E_0$  takes the form 
\[
E_0=\f12\displaystyle\int_{\Om_t}|v|^2dX+g\displaystyle\int_{\Om_t}X\cdot e_z dX+\sigma S(u)+[\gamma]X\cdot \tau_b|_{X_c}
\]
where $S(u)$ is the surface potential defined in \eqref{S(u) def}, and $[\gamma]u\cdot \tau_b|_{Y_c}$ is the interaction energy at the corner.

A standard way to analyze  the problem $\mbox{(WW)}$ is  linearization. We can start with the basic energy and try to find out the leading-order operator. Since it turns out  that the velocity part and the gravity part in $E_0$ are lower-order terms (see \cite{SZ, SZ2},  and also verified in \cite{Poyferre}), we focus on the last two terms. 

In fact, we denote
\[
E_c=\sigma S(u)+[\gamma]X\cdot \tau_b|_{X_c}
\] 
as the part of the energy related to the leading-order operator as well as  the contact point. We will compute variations of $E_c$.

Firstly, given a path $u(s,t,\cdot)\in \G$ and $\bar w=u_s|_{s=0}\in T_{u}\G$ with $w=\bar w\circ u^{-1}|_{s=0}$, we compute the  first variation of $E_c$  as  
\beq\label{1st variation of E0}
\langle E'_c,\bar w\rangle =\sigma\langle S'(u),\bar w\rangle +[\gamma]w\cdot \tau_b|_{X_c}
\eeq
where we already know from last subsection that
\[
\langle S'(u),\bar w\rangle=\displaystyle\int_{\G_t}\kappa w^\perp ds+w\cdot \tau_t|_{X_c}.
\]

Secondly, for the second variation, we also start with the surface potential $S(u)$. Let $h(s,\cdot)$ be  a geodesic on $\G$ with $h(0)=u$ and $\bar w=h_s|_{s=0}$, which means 
\[
\bar \cD_s\bar w=(D_s w+\na P_{w,w})\circ h=0.
\] 
Based on the first variation, the second variation for $S(u)$ can be written as 
\[
\bar \cD^2 S(u)(\bar w,\bar w)=\f d{ds}\displaystyle\int_{\G_t}\kappa w^\perp ds\big|_{s=0}+\f d{ds}w\cdot\tau_t|_{X_c}\big|_{s=0}:=A+B.
\]
A direct computation leads to
\[
\begin{split}
A=&\displaystyle\int_{\G_t}(D_s\kappa\, w\cdot n_t+\kappa D_sw\cdot n_t+\kappa\,w\cdot D_s n_t)ds+\kappa\,w\cdot n_t D_s ds\big|_{s=0}\\
=&\displaystyle\int_{\G_t}\Big(\big(-\D_{\G_t}w^\perp-w^\perp|\Pi|^2+\cD\cdot\Pi(w^\top)\big)w^\perp-\kappa\na P_{w,w}\cdot n_t-\kappa\na_{w^\top}w\cdot n_t\\
&\quad+\kappa w^\perp(\kappa w^\perp+\cD\cdot w^\top)\Big)ds
\end{split}
\] 
where one applied  equation \eqref{Dt k 1} for $D_s \kappa$ and 
\[
D_sw=-\na P_{w,w},\quad D_sn_t=-\big((\na w)^*n_t\big)^\top,\quad D_s ds=(\kappa w^\perp+\cD\cdot w^\top)ds.
\] 
Checking carefully on the term $\int_{\G_t}\big(\cD\cdot\Pi(w^\top)\big)w^\perp \,ds$ by integrating by parts and on $P_{w,w}$ term from Proposition \ref{prop:P(f,v)}, one can tell that the leading-order term is the first one in $A$. So one writes
\[
A=\displaystyle\int_{\G_t}(-\D_{\G_t}w^\perp)w^\perp ds\,+\,\hbox{lower-order terms.}
\] 
Moreover, in 2 dimensional case we can write that
\[
\begin{split}
\displaystyle\int_{\G_t}(-\D_{\G_t}w^\perp)w^\perp ds=&
-\displaystyle\int_{\G_t}\cD\cdot(\na^\top w^\perp)\,w^\perp ds\\
=&\displaystyle\int^{+\infty}_0\f d{ds}(\na_{\tau_t}w^\perp)\,w^\perp ds\\
=&-w^\perp\,\na_{\tau_t}w^\perp\big|_{X_c}+\displaystyle\int_{\G_t}|\na_{\tau_t}w^\perp|^2ds
\end{split}
\] 
which implies 
\[
A=\displaystyle\int_{\G_t}|\na_{\tau_t}w^\perp|^2ds-w^\perp\,\na_{\tau_t}w^\perp\big|_{X_c}\,+\,\hbox{lower-order terms.}
\]
Next, we turn to deal with term  $B$ to find 
\[
\begin{split}
B=&D_s w\cdot \tau_t+w\cdot D_s\tau_t\big|_{X_c}=-\na P_{w,w}\cdot \tau_t+w\cdot(\na_{\tau_t}w\cdot n_t)n_t\big|_{X_c}\\
=&-\na P_{w,w}\cdot \tau_t+w^\perp\,\na_{\tau_t}w^\perp-w^\perp(w\cdot\na_{\tau_t}n_t)\big|_{X_c},
\end{split}
\] 
which infers that
\[
B=w^\perp\,\na_{\tau_t}w^\perp|_{X_c}\,+\,\hbox{lower-order terms}
\]
by noticing $\na P_{w,w}$ is also a lower-order part from Proposition \ref{prop:P(f,v)}.

As a result, combining the analysis on $A,B$ we can conclude that 
\[
\bar \cD^2S(u)(\bar w,\bar w)=\displaystyle\int_{\G_t}|\na_{\tau_t}w^\perp|^2ds\,+\,\hbox{lower-order terms.}
\]

On the other hand, we consider the second variation for the second term in $E_c$ to find that
\[
\f{d}{ds}\big([\gamma] w\cdot \tau_b|_{X_c}\big)\big|_{s=0}=[\gamma]\big(D_s w\cdot \tau_b+w\cdot D_s \tau_b\big)\big|_{X_c}
=-[\gamma]\na P_{w,w}\cdot \tau_b|_{X_c},
\] where $D_s \tau_b=0$ since $\tau_b$ is constant near the contact point. Consequently,  one can see that this term  also turns out to be a lower-order term.

Summing up the analysis above, we finally  arrive at  the second variation for $E_c$ as  
\[
\bar\cD^2 E_c(u)(\bar w,\bar w)=\displaystyle\int_{\G_t}|\na_{\tau_t}w^\perp|^2ds\,+\,\hbox{lower-order terms.}
\] 
Consequently, we can see that,  in the presence of the contact point, the leading-order part from the second variation of our basic energy $E_0$ remains the same as the related classical case (see \cite{SZ}), while  the terms related to the corner are merely lower-order ones. 

This computation  tells us that in the linearization of system $\mbox{(WW)}$,  the leading-order operator should remain the same as the classical case,  which will be verified in the following text. Moreover, one can see in our paper that,  higher-order terms related to the corner appear only in the dissipation part.

\bigskip

\section{Dissipation equation}
\setcounter{equation}{0}

In this section, we will show that  system $\mbox{(WW)}$ satisfies the following  energy-dissipation equality 
\beq\label{dissipation eqn}
\f d{dt} E_0+\beta_c|v_c|^2=0.
\eeq
Compared to classical water waves problems, our energy has a dissipation related to the contact point, which is completely new. Moreover, dissipations involving the contact point also take place in our higher-order energy in the following sections.  

To prove the dissipation equality, recalling \eqref{basic energy} for the basic energy $E_0$,  we firstly compute  that
\[
\f{d}{dt} E_0=\displaystyle\int_{\Om_t} v\cdot D_t v dX+g\displaystyle\int_{\Om_t}D_t zdX+\sigma\p_tS(u)+[\gamma]v\cdot \tau_b|_{X_c},
\]
where 
\[
\p_tS(u)=\displaystyle\int_{\Om_t}\na\kappa_{\cH}\cdot vdX+v\cdot \tau_t|_{X_c}.
\]
Noticing that 
\[
v_c=-v\cdot \tau_b|_{X_c},\quad v\cdot\tau_t|_{X_c}=v_c\cos \om,
\]
we arrive at 
\[
\f{d}{dt} E_0=\displaystyle\int_{\Om_t}v\cdot(D_t v+g e_z+\sigma\na \kappa_{\cH})dX
+\big(\sigma v_c\cos\om-[\gamma] v_c\big)\big|_{X_c}.
\]
Plugging in the Euler equation and condition \eqref{condition at X_c} from $\mbox{(WW)}$, we find that
\[
\begin{split}
\f{d}{dt} E_0=&-\displaystyle\int_{\Om_t}v\cdot\na P_{v,v}dX+ \sigma v_c\big(\cos\om-\f{[\gamma]}{\sigma}\big)\big|_{X_c}\\
=&-\beta_c|v_c|^2\big|_{X_c}
\end{split}
\] 
and the proof ends.

\bigskip

\section{Elliptic estimates and trace theorems in the corner domain}
\setcounter{equation}{0}

In order to perform the energy estimates, some technical preparations are needed. This section provides us some useful elliptic estimates and trace theorems on the corner domain, which are  adjusted from \cite{PG1,MW}.  

To begin with, the following mixed-boundary elliptic system are used frequently in our paper:
\beq\label{eq:elliptic}
\left\{
\begin{array}{l}
\Delta u=h,\qquad \textrm{on}\quad \Om_t, \\
u|_{\Gamma_t}=f,\quad\quad \na_{n_b} u+b_0\na_{\tau_b}  u|_{\Gamma_b}=g,
\end{array}
\right.
\eeq
where $b_0$ is a constant coefficient. 
\begin{theorem}\label{est:elliptic}
{\it If the contact angle $\om(t)\in(0,\pi/4)$, we have the following estimate for system \eqref{eq:elliptic} when $2\leq s\leq 3$:
\beno
\|u\|_{H^{s}(\Om_t)}\leq C(\| \Gamma_t\|_{H^{\f52} \cap H^{s-\f12}  },\| \Gamma_b\|_{H^{\f52}    })(\|h\|_{H^{s-2}(\Om_t)}+ \|f\|_{H^{s-\f12}(\Gamma_t)}+\|g \|_{H^{s-\f32}(\Gamma_b)}).
\eeno
  For $s=4$ and $\om(t)\in(0,\pi/6)$, we have
\beno
\|u\|_{H^{s}(\Om_t)}\leq C(\| \Gamma_t\|_{ H^{\f72}  },\| \Gamma_b\|_{  H^{\f72}   })(\|h\|_{H^{2}(\Om_t)}+ \|f\|_{H^{\f72}(\Gamma_t)}+\|g \|_{H^{\f52}(\Gamma_b)}).
\eeno
} 
\end{theorem}
 \no {\bf Proof.}
When $s=2, 3, 4$, The estimate comes directly from Proposition 5.19 and Remark 5.20 \cite{MW}. For the case $2<s<3$,  a complex interpolation is applied to finish the proof.

\ef

Meanwhile, we also give an a  priori estimate for  \eqref{eq:elliptic} when $b_0$ is a  function with a bounded support on $\G_b$, which will be used in Proposition \ref{prop:v}. 
\begin{lemma}\label{est:elliptic-1}
{\it Let  the contact angle $\om(t) \in(0,\pi/4)$ and $u\in H^3(\Om_t)$ be  a solution of system \eqref{eq:elliptic} with $h\in H^1(\Om_t)$, $f\in H^\f52(\G_t)$ and $g\in H^\f32(\G_b)$. Moreover, assume  that  the  coefficient $b_0\in H^\f32(\G_b)$ is a function vanishing away from the contact point. Then one has the following elliptic estimate
\[
\|u\|_{H^{3}(\Om_t)}\leq C(\| \Gamma_t\|_{H^{\f52}   },\| \Gamma_b\|_{H^{\f52}   }, \|b_0\|_{H^{\f32}(\Gamma_b)})(\|h\|_{H^{1}(\Om_t)}+ \|f\|_{H^{ \f52}(\Gamma_t)}+\|g \|_{H^{ \f32}(\Gamma_b)}).
\]
}
\end{lemma}

\no {\bf Proof.}
The proof is based on Theorem \ref{est:elliptic} and we will use a unit decomposition to localize system 
\eqref{eq:elliptic}.  The constant $C$ in the proof stands for $C(\| \Gamma_t\|_{H^{\f52}   },\| \Gamma_b\|_{H^{\f52}   }, \|b_0\|_{H^{\f32}(\Gamma_b)})$.

To begin with,  since $b_0$ has a bounded support, we choose a unit decomposition $\sum^k_{i=1}\chi_i=1$ with some $k\in \N$ such that $\chi_1$ is supported near $X_c$ and $b_0$ vanishes on the support of $\chi_k$. Moreover, the horizontal size of the supports for $\chi_1,\dots,\chi_{k-1}$ is a small constant $\delta>0$ to be fixed later. As a result, based on the decomposition
\[
u=\sum^k_{i=1}\chi_i u:=\sum^k_{i=1}u_i,
\] 
we need to prove the desired estimate for each $u_i$ to close the proof. 
In fact, it suffices to  focus on the estimate for $u_1$ since the remaining parts will be similar and classical. 

A direct computation leads to the system for $u_1$:
\[
\left\{
\begin{array}{l}
\Delta u_1=h+[\Delta,\chi_1]u,\quad \textrm{on}\quad \Om_t, \\
u_1|_{\Gamma_t}=\chi_1 f,\\
\na_{n_b} u_1+b_0(X_c)\na_{\tau_b} u_1\big|_{\Gamma_b}= g_1
\end{array}
\right.
\]
where 
\[
g_1=\chi_1   g+ u\na_{n_b} \chi_1  +b_0\,u\na_{\tau_b} \chi_1  +(b_0(X_c)-b_0)\na_{\tau_b} u_1\big|_{\Gamma_b}.
\]
In order to prove the estimate for $u_1$, we firstly deal with $g_1$ on  $\Gamma_b$ to arrive at
\[
\|g_1\|_{H^\f32(\G_b) }\leq C \|g \|_{H^{ \f32}(\Gamma_b)}+C (1+\|b_0\|_{H^{\f32}(\Gamma_b)}) \|u\|_{H^{2}(\Om_t)} +C\|b_0(X_c)-b_0\|_{L^\infty(supp\chi_1)} \|u_1\|_{H^3(\Om_t)}
\]
where Theorem \ref{trace thm PG}  is used on $u$.
Applying Theorem \ref{est:elliptic},  we have 
\[
\begin{split}
\|u_1\|_{H^3(\Om_t)} \le &C(\|h\|_{H^{1}(\Om_t)}+\|u\|_{H^{2}(\Om_t)}+ \|f\|_{H^{ \f52}(\Gamma_t)}+\|g \|_{H^{ \f32}(\Gamma_b)})\\
&\quad+\|b_0(X_c)-b_0\|_{L^\infty(supp \chi_1)} \|u_1\|_{H^3(\Om_t)}. 
\end{split}
\]
Consequently, when the horizontal size $\delta$ of $supp\chi_1$ is small enough,   a small constant $C_\delta$ can be found such that
\beno
\|b_0(X_c)-b_0\|_{L^\infty(supp \chi)}\leq C_\delta,
\eeno
which leads to the  estimate for $u_1$:
\[
\|u_1\|_{H^3(\Om_t)}\le C\big(\|h\|_{H^{1}(\Om_t)}+\|u\|_{H^{2}(\Om_t)}+ \|f\|_{H^{ \f52}(\Gamma_t)}+\|g \|_{H^{ \f32}(\Gamma_b)}\big). 
\]
Moreover,  $\|u\|_{H^2(\Om_t)}$ on the right side of the above inequality can be handled by an interpolation, Theorem \ref{trace thm PG} and the Poincar\'e inequality.  Therefore,  our proof is finished.
\ef

Except for the elliptic estimates above, we present here some useful  trace theorems in the corner domain.
\begin{theorem}\label{trace thm PG} (Traces on $\G_t$ or $\G_b$)
{\it Let the integer $ l\in[0,s-\f12)$ with some $s>l+\f12$, we define the map
\beno
u\to \{u, \na_{n_j} u,\dots \na_{n_j}^l u\}|_{\Gamma_j},
\eeno
 for $u\in \mathscr D(\bar \Om_t)$ where $n_j$ is the unit outward normal vector on $\Gamma_j$ with $\Gamma_j$ taking $\Gamma_b$ or $\Gamma_t$. Then, the map has a unique continuous extension as an operator from
\beno
H^{s}(\Om_t)\quad \textrm{onto} \quad  \Pi_{i=0}^lH^{s-i-\f12}(\Gamma_j).
\eeno
Moreover, one has the estimate for $0\le i\le l$:
\[
\| \na^i_{n_b} u\|_{H^{s-i-\f12 }(\Gamma_b )}+\|  \na^i_{n_t}u\|_{H^{s-i-\f12 }(\Gamma_t )}\le    C(\|\Gamma_t\|_{H^{s-\f12} })\|u\|_{H^{s}(\Om_t)}.
\]
}
\end{theorem}
\no{\bf Proof.} This result is adjusted from Theorem 1.5.2.1 \cite{PG1} and Remark 4.2 \cite{MW} by a  cut-off function argument and interpolations.
\ef

\begin{theorem}\label{thm:trace} (Trace theorem with mixed boundary conditions) 
{\it  Let $m\in [2, 4]$ be an integer and  functions $f\in H^{m-\f12}(\Gamma_t)$, $g \in H^{m-\f 32}(\Gamma_b)$ be given. Then there exists a function $u\in H^{m}(\Om_t)$ satisfying the following mixed boundary conditions
\beno
u|_{\Gamma_t}=f,\quad \na_{n_b} u+b_0\na_{\tau_b}u|_{\Gamma_b}=g 
\eeno
where $b_0$ is a constant coefficient and  ${\bf n}_b+b_0{\bf \tau}_b\nparallel{\bf\tau}_t$.
Moreover, one has the estimate
\ben\label{111}
\|u\|_{H^m(\Om_t )}\le    C(\|\Gamma_t\|_{H^{m-\f12} }, \|\Gamma_b\|_{H^{m-\f12} })(\|f \|_{H^{m-\f12}(\Gamma_t)}+\|g \|_{H^{m-\f32}(\Gamma_b)})
\een
and
\ben\label{222}
\|u\|_{H^{m-\f12}(\Om_t )}\le    C(\|\Gamma_t\|_{H^{m-\f12} }, \|\Gamma_b\|_{H^{m-\f12} })(\|f \|_{H^{m-1}(\Gamma_t)}+\|g \|_{H^{m-2}(\Gamma_b)}).
\een
}
\end{theorem}
\no {\bf Proof.}
The estimate \eqref{111} can be found  in Theorem 4.6 \cite{MW}.  For the estimate  \eqref{222}, it can be proved by a complex interpolation theorem.  Besides, the condition ${\bf n}_b+b_0{\bf \tau}_b\nparallel{\bf\tau}_t$ can always be satisfied in this paper.
\ef

Next,  some special trace theorems involving $\tilde H^\f12(\G_b)$ and $\tilde H^{-\f12}(\G_b)$ are also needed in our paper. 
\begin{lemma}\label{H1 zero trace}
{
\it  Assume that $u|_{\G_t}=0$ and let $f=u|_{\G_b}$ for any function $u\in H^1(\Om_t)$. Then the mapping $u\mapsto f$  is  linear continuous from $H^1(\Om_t)$ onto  $\tilde H^\f12(\G_b)$:
\[
\|f\|_{\tilde H^\f12(\G_b)}\le C(\|\Gamma_t\|_{H^{\f52} }) \|u\|_{H^1(\Om_t)}.
\]
}
\end{lemma}
\no {\bf Proof.}
This is an adaption from Theorem A.1 \cite{MW}.  In fact, using a cut-off function $\chi$ near the contact point,  one can prove the estimate for $\chi u$ by using Theorem A.1 \cite{MW}:
\[
\| \chi f \|_{\tilde H^\f12(\G_b)}\le C(\|\Gamma_t\|_{H^{\f52} })\|\chi u\|_{H^1(\Om_t)}
\]
where the dependence of the coefficient $C$ can be checked from Section 1.5.2 in \cite{PG1}.
 The remainder part $(1-\chi)u$ can be dealt with a classical trace theorem to have 
\[
\|(1-\chi) f \|_{H^\f12(\G_b)}\le C(\|\Gamma_t\|_{H^{\f52} } ) \|(1-\chi )u\|_{H^1(\Om_t)}.
\]
Combining these two estimates, the proof is finished. 
\ef

\begin{lemma}\label{G_b H-1/2}
{\it If $u$ belongs to 
$H^\f12(\G_b)$, then $\na_{\tau_b} u$ belongs to  $\tilde H^{-\f12}(\G_b)$ and satisfies the estimate
\[
\|\na_{\tau_b}u\|_{\tilde H^{-\f12}(\G_b)}\le C(\|\Gamma_t\|_{H^{\f52} } ) \|u\|_{H^\f12(\G_b)}.
\] 
}
\end{lemma}
\no{\bf Proof.} This lemma is an application of Theorem A.2 \cite{MW} and the proof of Theorem 1.4.4.6 \cite{PG1}, where the dependence of the coefficient $C$ can also be checked from the proofs of these theorems.
\ef

\begin{lemma}\label{na_n_b w on G_b}
{\it Let $u\in E(\D; L^2(\Om_t))=\{u\in H^1(\Om_t)|\,\D u\in L^2(\Om_t)\}$. Then the mapping
$u\mapsto \na_{n_b}u|_{\G_b}$
is continuous  from $E(\D;L^2(\Om_t))$ into $\tilde H^{-\f12}(\G_b)$:
\[
\| \na_{n_b}u \|_{\tilde H^{-\f12}(\G_b)}\le C(\|\Gamma_t\|_{H^{\f52} } ) \big(\|u\|_{H^1(\Om_t)}+\|\D u\|_{L^2(\Om_t)}\big).
\]
}
\end{lemma}
\no{\bf Proof.} It can be proved  directly from Theorem 1.5.3.10 \cite{PG1}, while one can check the coefficient $C$ as before.
\ef

It turns out that the Sobolev's embedding theorem also works on  the corner domain, and we only pick up three cases here.
\begin{lemma}\label{embedding}
{\it We have the following embeddings:
\[
\|u\|_{L^4(\Om_t)}\le C(\|\Gamma_t\|_{H^{\f52} })  \|u\|_{H^\f12(\Om_t)}
\] 
for any $u\in H^\f12(\Om_t)$,  and 
\[
\|u\|_{L^\infty(\Om_t)}\le C(\|\Gamma_t\|_{H^{\f52} }) \|u\|_{H^{s_1}(\Om_t)}
\] 
for  any $u\in H^{s_1}(\Om_t)$ with $s_1=1+\epsilon$ ($\epsilon>0$ is a small constant). Moreover, for any $ f\in H^{s_2}(\G_b)$ with $s_2=\f12+\epsilon$, the embedding holds:
\[
\| f \|_{L^\infty(\G_b)}\le C(\|\Gamma_t\|_{H^{\f52} }) \|f\|_{H^{s_2}(\G_b)}.
\]
}
\end{lemma}
\no{\bf Proof.}
Using an extension theorem (for example Theorem 1.4.3.1 \cite{PG1}) firstly to remove the corner and then the classical Sobolev's embedding theorem, this lemma can be proved.
\ef

Moreover, except for the trace theorems and  embedding results above, we will meet  with $H^1$-type elliptic estimates frequently in the energy estimates. The following lemma is about the $H^1$ estimate related to the harmonic extension operator $\cH$.  
\begin{lemma}\label{Harmonic extension H1 estimate}
{\it For a given function $f\in H^\f12(\G_t)$, the system 
\beq\label{H f system}
\left\{\begin{array}{ll}
\D \cH(f)=0,\qquad\hbox{on}\quad \Om_t,\\
\cH(f)|_{\G_t}=f,\quad \na_{n_b}\cH(f)|_{\G_b}=0.
\end{array}\right.
\eeq
admits a  solution  $\cH(f)\in H^1(\Om_t)$ as  the harmonic extension of $f$,  and the estimate below holds:
\[
\|\cH(f)\|_{H^1(\Om_t)}\le C(\|\Gamma_t\|_{H^{\f52} }) \|f\|_{H^\f12(\G_t)}.
\]
}
\end{lemma}
\no{\bf Proof.}
\no Step 1: Find $u_1\in H^1(\Om_t)$ satisfying
\[
\left\{\begin{array}{ll}
\D u_1=0,\qquad\hbox{on}\quad \Om_t,\\
u_1|_{\G_t}=f.
\end{array}\right.
\]
In fact, the only boundary condition here is on $\G_t$, so one can extend $\G_t$ and $f$ such that  $f^{ex}$ is defined on a  horizontally infinite curve $\G^{ex}_t$.  Moreover, one can also extend the domain $\Om_t$ to a horizontally infinite strip $\Om^{ex}_t$ as the domain in \cite{Lannes} and pose zero Dirichlet boundary condition on the bottom if necessary. Therefore, the system for $u_1$ can be extended to a system for $u^{ex}_1$ on $\Om^{ex}_t$. Applying Lemma 2.12 \cite{Lannes} leads to the estimate 
\[
\|u^{ex}_1\|_{H^1(\Om^{ex}_t)}\le C( \|\Gamma_t\|_{H^{\f52} }) \,\|f^{ex}\|_{H^\f12(\G^{ex}_t)},
\]
which infers that $u_1=u^{ex}_1|_{\Om_t}$ solves the system above and 
\[
\|u_1\|_{H^1(\Om_t)}\le C(\|\Gamma_t\|_{H^{\f52} } ) \,\|f\|_{H^\f12(\G_t)}.
\]

\no Step 2: Define $u_2$ by the system 
\[
\left\{\begin{array}{ll}
\D u_2=0,\qquad\hbox{on}\quad \Om_t,\\
u_2|_{\G_t}=0,\quad \na_{n_b}u_2|_{\G_b}=-\na_{n_b}u_1|_{\G_b}.
\end{array}\right.
\]
According to Lemma \ref{na_n_b w on G_b}, we  know that $\na_{n_b}u_1|_{\G_b}\in \tilde H^{-\f12}(\G_b)$, so the boundary conditions make sense. The equivalent variation equation for $u_2$ is
\[
\displaystyle\int_{\Om_t}\na u_2\cdot \na \phi dX=-\displaystyle\int_{\G_b}\na_{n_b}u_1\,\phi ds
\]
with any $\phi\in V=\{\phi\in H^1(\Om_t)\,\big|\,\phi|_{\G_t}=0\}$. From Lemma \ref{H1 zero trace}, we know that $\phi\in \tilde H^{\f12}(\G_b)$, so the right-hand side of the variation equation satisfies
\[
\big|\displaystyle\int_{\G_b}\na_{n_b}u_1\,\phi ds\big|
\le \|\na_{n_b}u_1\|_{\tilde H^{-\f12}(\G_b)}\|\phi\|_{\tilde H^{\f12}(\G_b)}
\le C(\|\Gamma_t\|_{H^{\f52} }) \, \|u_1\|_{H^1(\Om_t)}\|\phi\|_{H^1(\Om_t)}.
\]
As a result, From Lax-Milgram theorem there exists a unique solution $u_2$ solving the variational system with the estimate
\[
\|u_2\|_{H^1(\Om_t)}\le C(\|\Gamma_t\|_{H^{\f52} })\,\|u_1\|_{H^1(\Om_t)}\le C(\|\Gamma_t\|_{H^{\f52} })\,\|f\|_{H^\f12(\G_t)},
\]
where  the Poincar\'e inequality is applied. 

\no Step 3: Summing the first two steps together and  let
\[
\cH(f)=u_1+u_2,
\]
we conclude that $\cH(f)\in H^1(\Om)$ solves system \eqref{H f system} and the desired estimate holds.
\ef
   
 In the end of this section, we consider elliptic estimates for a Neumann-type system  
 \ben\label{eq:elliptic-N}
\left\{
\begin{array}{l}
\Delta u=h,\quad \textrm{on}\quad \Om_t, \\
\na_{n_t}u|_{\Gamma_t}=f,\quad\quad \na_{n_b} u |_{\Gamma_b}=g
\end{array}
\right.
\een
satisfying the compatibility condition
\beno
\int_{\Om_t} hdX=\int_{\Gamma_t} fds+\int_{\Gamma_b}gds.
\eeno
One can check from \cite{MW} with the Neumann conditions above and conclude the following result for system \eqref{eq:elliptic-N}. The proof is similar as that for Proposition 5.19 \cite{MW} without singularity and hence is omitted. 
\begin{theorem}\label{est:elliptic-N}
{\it 
Let the contact angle $\om(t)\in(0,\pi/4)$, $h\in H^{\f12}(\Om_t)$, $f\in H^{1}(\Gamma_t)$ and $g\in H^{1}(\Gamma_b)$. If $u\in H^\f52(\Om_t)$ is the solution to system \eqref{eq:elliptic-N},   the following estimate holds:
\beno
\|u\|_{H^{\f52}(\Om_t)}\leq C(\| \Gamma_t\|_{H^{\f52}})(\|h\|_{H^{\f12}(\Om_t)}+ \|f\|_{H^{1}(\Gamma_t)}+\|g \|_{H^{1}(\Gamma_b)}).
\eeno
}
\end{theorem}

\bigskip
  
\section{New formulation for $J$}
\setcounter{equation}{0}
  
In this section, we will derive the equation for a good unknown $J=\na \kappa_{\cH}$ which is introduced in \cite{SZ}. Recalling that
\[
\na P=J+\na P_{v,v},
\] we can see that $J$ is the part related to the mean curvature $\kappa$ in the pressure term $\na P$, which turns out to be the main part compared to $\na P_{v,v}$.

To begin with, we need some expansions of some commutators. The computations here follow the formulation of Shatah-Zeng, see \cite{SZ} (some expressions are quoted from there directly). In order to be self-contained,  we will also recall some details from their work.

\subsection{Commutators and  computations about $\kappa$} 
\noindent 1. $[D_t,\cH]$.  Firstly, recall that  $\D^{-1}(h,g)$ is defined as  the solution $u$ to the system
\beq\label{del-1 system}
\left\{\begin{array}{ll}
\D u=h\qquad \hbox{on}\quad \Om_t\\
u|_{\G_t}=0,\qquad \na_{n_b}u|_{\G_b}=g.
\end{array}
\right.
\eeq

In order to analyze the commutator, we start with the elliptic system of $D_tf_{\cH}$ where recall also that $f_{\cH}$ is the harmonic extension  of $f$ defined $\G_t$. In fact, $D_tf_{\cH}$ satisfies the equation 
\[
\D D_t f_{\cH}=[\D,D_t]f_{\cH}
\]
and the boundary condition on $\G_b$:
\[
\na_{n_b}D_tf_{\cH}|_{\G_b}=[\na_{n_b},D_t]f_{\cH}|_{\G_b}.
\]
In a word, direct computations shows that $D_t f_{\cH}$ solves the system
\[
\left\{\begin{array}{ll}
\D D_tf_{\cH}=2\na v\cdot \na^2 f_{\cH}+\D v\cdot \na f_{\cH}\qquad \hbox{on}\quad \Om_t\\
D_tf_{\cH}|_{\G_t}=D_t f,\qquad \na_{n_b}D_tf_{\cH}|_{\G_b}=(\na_{n_b}v-\na_v n_b)\cdot\na f_{\cH}|_{\G_b}.
\end{array}
\right.
\] 
As a result,  we can have 
\beq\label{commutator Dt H}
D_t f_{\cH}=\cH( D_t f)+\D^{-1}\big(2\na v\cdot\na^2f_{\cH}+\D v\cdot\na f_{\cH},\, (\na_{n_b}v-\na_v n_b)\cdot \na f_{\cH}\big)
\eeq

\noindent 2. $[D_t,\,\D^{-1}]$. Denoting $u=\D^{-1}(h,g)$,  we know that $u$ solves  system \eqref{del-1 system}. To compute this commutator, we investigate the system for $D_t u$.

In fact, direct computations shows that $D_t u$ satisfies 
 \[
\left\{\begin{array}{ll}
\D D_tu=D_t h+2\na v\cdot \na^2 u+\D v\cdot \na u\qquad \hbox{on}\quad \Om_t\\
D_tu|_{\G_t}=0,\qquad \na_{n_b}D_tu\big|_{\G_b}=D_t g+(\na_{n_b}v-\na_v n_b)\cdot\na u\big|_{\G_b}.
\end{array}
\right.
\] 
Consequently,  we conclude by  a direct decomposition that
\beq\label{commutator Dt D-1}
\begin{split}
D_t\D^{-1}(h,g)=&\D^{-1}(D_th,\,D_tg)\\
&\,
+\D^{-1}\Big(2\na v\cdot\na^2\D^{-1}(h,g)+\D v\cdot \na \D^{-1}(h,g),\,
(\na_{n_b}v-\na_{v}n_b)\cdot \na\D^{-1}(h,g)\Big).
\end{split}
\eeq

\noindent 3. $[D_t,\,\cN]$. The Dirichlet-Neumann operator $\cN$ is defined by
\[
\cN f:=\na_{n_t}f_{\cH}\qquad\hbox{on}\quad \G_t
\] 
for a given function $f$ on $\G_t$. We know directly from Theorem 1.2 \cite{MW} that
\beq\label{commutator DN}
\begin{split}
[D_t,\,\cN]f=&\na_{n_t}\D^{-1}\Big(2\na v\cdot \na^2f_{\cH}+\D v\cdot\na f_{\cH},\,(\na_{n_b}v-\na_{v}n_b)\cdot\na f_{\cH}\Big)\\
&\,-\na_{n_t}v\cdot \na f_{\cH}-\na_{(\na f_{\cH})^\top} v\cdot n_t\qquad\qquad\hbox{on}\quad \G_t.
\end{split}
\eeq 

\noindent   4. $[D_t,\,\D_{\G_t}]$.  One can have by a direct computation that
\beq\label{commutator surface delta}
D_t\D_{\G_t}f=\D_{\G_t}D_t f+2\cD^2f\big(\tau_t,(\na_{\tau_t}v)^\top\big)-(\na f)^\top\cdot \D_{\G_t}v+\kappa \na_{(\na f)^\top} v\cdot n_t\qquad\hbox{on}\quad \G_t.
\eeq

\noindent 5. $D_t\kappa$. Direct computations lead to the following expression
\beq\label{Dt k 1}
D_t\kappa=-\D_{\G_t}v^\perp-v^\perp|\Pi|^2+\cD\cdot \Pi(v^\top)
\eeq
or equivalently
\beq\label{Dt k 2}
D_t\kappa=-\D_{\G_t}v\cdot n_t-2\Pi(\tau_t)\cdot \na_{\tau_t}v \qquad\hbox{on}\quad \G_t.
\eeq

\noindent 6. $D_t^2\kappa$. Based on $D_t\kappa$, we  can have
\beq\label{Dt2 kappa}
D^2_t\kappa=-n_t\cdot \D_{\G_t}D_tv+2\s \Pi(\tau_t)\cdot \na_{\tau_t} J+R_1\qquad\hbox{on}\quad \G_t
\eeq   with 
\[
\begin{split}
R_1=& 2\Big[\cD\big((\na v)^*n_t\big)^\top+\Pi\big((\na_{\tau_t} v)^\top\big)\Big]\cdot \na_{\tau_t} v+\D_{\G_t}v\cdot\big((\na v)^*n_t\big)^\top+2\Pi(\tau_t)\cdot (\na v)^2\tau_t\\
&
\  -2(\na_{\tau_t} v\cdot n_t)\Big(\Pi(\tau_t)\cdot\na_{n_t}v
+\Pi(\tau_t)\cdot(\na v)^*n_t\Big) 
-2n_t\cdot D^2v\big(\tau_t,\,(\na_{\tau_t}v)^\top\big)\\
&
\ -2(\na_{n_t}v\cdot n_t)(\Pi(\tau_t)\cdot \na_{\tau_t} v)+n_t\cdot \na v\big((\D_{\G_t}v)^\top\big) -\kappa\big|\big((\na v)^*n_t\big)^\top\big|^2+2\Pi(\tau_t)\cdot \na_{\tau_t}\na P_{v,v}
\end{split}
\] 
where one can see that the leading-order terms in $R_1$ are like $\na^2v,\na n_t$.

\subsection{Equations for $J$}

With the preparations above, we are ready to derive  equations for $J$.  In the following analysis,  the idea is trying to express $D_tJ, D^2_t J$  using $v$, $D_tv$, $J$ and $D_t J$ instead of $\kappa$, $D_t\kappa$ or $D^2_t\kappa$.

To get started, one recalls that $J=\na\kappa_{\cH}$ and applies $D_t$  to find
\[
\begin{split}
D_t J=&D_t\na\kappa_{\cH}=\na D_t\kappa_{\cH}-(\na v)^*J\\
=&\na\cH(D_t\kappa)+\na[D_t,\cH]\kappa-(\na v)^*J
\end{split}
\]
where we used the commutator
\[
[D_t,\na]=-(\na v)^*\na. 
\]
Plugging \eqref{Dt k 2} and \eqref{commutator Dt H} into the equation above,  we can find the equation for $D_tJ$:
\beq\label{Dt J}
\begin{split}
D_t J=&\na\cH(D_t\kappa)+\na\D^{-1}\big(2\na v\cdot \na J
+\D v\cdot J,\,(\na_{n_b}v-\na_v n_b)\cdot J\big)-(\na v)^*J\\
=&-\na \cH\big(\D_{\G_t}v\cdot n_t+2\Pi(\tau_t)\cdot \na_{\tau_t}  v\big)+\na\D^{-1}\big(2\na v\cdot \na J
+\D v\cdot J,\,(\na_{n_b}v-\na_v n_b)\cdot J\big)\\
&\quad-(\na v)^*J.
\end{split}
\eeq

Secondly, to derive the equation for $D^2_tJ$, we begin with
\[
\begin{split}
D^2_tJ=&D^2_t\na\kappa_{\cH}=D_t\big(\na D_t\kappa_{\cH}+[D_t,\na]\kappa_{\cH}\big)\\
=&\na D^2_t\kappa_{\cH}-(\na v)^*\na D_t\kappa_{\cH}-(D_t\na v)^*J-(\na v)^*D_t J\\
=&\na D^2_t \kappa_{\cH}-2(\na v)^*D_tJ-(\na D_t v)^* J-\big((\na v)^2\big)^*J+(\na v)^*\na vJ
\end{split}
\] 
where  
\[
\na D_t\kappa_{\cH}=D_t J+(\na v)^*J
\]
and 
\[
(D_t\na v)^*J=(\na D_t v)^*J-\big((\na v)^2\big)^*J.
\] 
 Now we need to deal with $D^2_t\kappa_{\cH}$ and relate it to $D^2_t\kappa$. In fact, direct computations lead to
\[
\begin{split}
D^2_t\kappa_{\cH}=& D_t\big(\cH(D_t\kappa)+[D_t,\cH]\kappa\big)\\
=&D_t\Big(\cH(D_t\kappa)+\D^{-1}\big(2\na v\cdot\na J+\D v\cdot J,(\na_{n_b}v-\na_vn_b)\cdot J\big)\Big)\\
=&\cH(D^2_t\kappa)+[D_t,\cH]D_t\kappa+D_t\D^{-1}\big(2\na v\cdot\na J+\D v\cdot J,(\na_{n_b}v-\na_vn_b)\cdot J\big)
\end{split}
\] 
where we applied \eqref{commutator Dt H} repeatedly. Plugging this expression above into the equation for $D^2_t J$, we arrive at 
\beq\label{Dt2 J eqn1}
D^2_tJ=\na\cH(D^2_t\kappa)+A_1+A_2+A_3
\eeq 
where
\[
A_1=\na [D_t,\cH]D_t\kappa,\quad A_2=\na D_t\D^{-1}\big(2\na v\cdot\na J+\D v\cdot J,(\na_{n_b}v-\na_vn_b)\cdot J\big)
\]
and 
\[
A_3=-2(\na v)^*D_tJ-(\na D_t v)^* J-\big((\na v)^2\big)^*J+(\na v)^*\na vJ.
\]
Applying \eqref{commutator Dt H}  and \eqref{Dt J} we find that
\beq\label{A_1}
\begin{split}
A_1=\na w=&\na\D^{-1}\big(2\na v\cdot\na+\D v\cdot ,(\na_{n_b}v-\na_vn_b)\cdot\big)\\
&\quad
\Big(D_tJ-\na\D^{-1}\big(2\na v\cdot \na J
+\D v\cdot J,\,(\na_{n_b}v-\na_v n_b)\cdot J\big)+(\na v)^*J\Big)
\end{split}
\eeq
where 
\[
w= \D^{-1}\big(2\na v\cdot\na^2\cH(D_t\kappa)+\D v\cdot\na \cH(D_t\kappa) ,(\na_{n_b}v-\na_vn_b)\cdot\na \cH(D_t\kappa)\big).
\]

On the other hand, applying \eqref{commutator Dt D-1} to $A_2$ we have 
\beq\label{A_2}
A_2=\na\D^{-1}\big(2\na v\cdot\na^2w_{A2}+\D v\cdot \na w_{A2},(\na_{n_b}v-\na_vn_b)\cdot\na w_{A2}\big) +\na\D^{-1}(h_{A2},\, g_{A2})
\eeq
where
\[
\begin{split}
 w_{A2}=&\D^{-1}\big(2\na v\cdot \na J+\D v\cdot J,\,(\na_{n_b}v-\na_v n_b)\cdot J\big),\\
 h_{A2}=&2\na v\cdot(\na D_tJ-(\na v)^*J)+2(\na D_tv-(\na v)^*\na v)\cdot \na J+D_tJ\cdot\D v\\
 &\quad+J\cdot (\D D_tv-\D v\cdot\na v-2\na v\cdot \na^2 v),\\
g_{A2}=&(\na_vn_b-\na_{n_b}v)\cdot \na v\cdot J+\na_{n_b}D_t v\cdot J-(D_t v-\na_v v)\cdot \na n_b\cdot J\\
&\quad +\na_v\big((\na v)^*n_b\big)^\top\cdot J+(\na_{n_b}v-\na_vn_b)\cdot D_tJ.
\end{split}
\]
One can tell that the leading-order terms in $A_1, A_2$  are like $J,D_tJ,\na v,\na D_tv$.

Now it remains to rewrite $\na \cH(D^2_t\kappa)$ in \eqref{Dt2 J eqn1}. In fact, substituting \eqref{Dt2 kappa} into this term, one has
\[
\begin{split}
\na\cH(D^2_t \kappa)=&\na\cH\Big(-n_t\cdot\D_{\G_t}D_tv+2\sigma\Pi(\tau_t)\cdot \na_{\tau_t} J+R_1\Big)\\
=&\na\cH\Big(\sigma n_t\cdot\D_{\G_t}J+ n_t\cdot\D_{\G_t}\na P_{v,v}+2\sigma\Pi(\tau_t)\cdot \na_{\tau_t}  J+R_1\Big)
\end{split}
\] 
where the Euler equation is applied. 
Moreover, direct computations show that
\[
2\sigma\Pi(\tau_t)\cdot \na_{\tau_t}  J=2\sigma\na_{\tau_t}n_t\cdot \na_{\tau_t}J=\sigma\D_{\G_t}J^\perp-\sigma\D_{\G_t}n_t\cdot J-\sigma n_t\cdot\D_{\G_t}J,
\] 
which leads to
\[
\na\cH(D^2_t\kappa)=\sigma \na\cH(\D_{\G_t}J^\perp)-\sigma\na\cH(\D_{\G_t}n_t\cdot J)+\na\cH\big(n_t\cdot\D_{\G_t}\na P_{v,v}\big)+\na\cH(R_1).
\]
Therefore, plugging the expression above into \eqref{Dt2 J eqn1}, we finally arrive at 
\beq\label{eqn for Dt2 J}
D^2_tJ=\s\na \cH(\D_{\G_t}J^\perp)+R_0.
\eeq 
with 
\[
R_0=-\s\na \cH(J\cdot \D_{\G_t}n_t)+\na \cH(n_t\cdot \D_{\G_t}\na P_{v,v})+\na\cH(R_1) +A_1+A_2+A_3
\]

\medskip
Based on the equations for $D_t J$, we will  derive one more equation for $\cD_tJ$ needed in the energy estimate. In fact, we know from the Hodge decomposition that 
\[
\cD_t J=D_t J+\na P_{J,v}
\] where $P_{J,v}$ satisfies the system 
\[\left\{\begin{array}{ll}
\D P_{J,v}=-tr (\na J\na v),\qquad\hbox{on}\quad\Om_t\\
P_{J,v}|_{\G_t}=0,\qquad \na_{n_b} P_{J,v}|_{\G_b}=\na_v n_b\cdot J.
\end{array}\right.
\] 
This implies 
\[
\begin{split}
D_t\cD_tJ= D_t(D_tJ+\na P_{J,v})=D^2_tJ+D_t\na P_{J,v}.
\end{split}
\]
Substituting \eqref{eqn for Dt2 J} into the equation above, we conclude the following equitation for $\cD_tJ$:
\beq\label{eqn cDt2 J}
D_t\cD_t J+\s \cA J=R
\eeq where the operator $\cA$ is defined by 
\[
\cA(w)=\na\cH\big(-\D_{\G_t}(w|_{\G_t})^\perp\big)
\]
and the remainder term $R$ is 
\[
R=R_0+D_t\na P_{J,v} .
\]
To close this section,  we  consider about $D_t\na P_{J,v}$ in $R$. 
Indeed, direct computations and \eqref{commutator Dt D-1} lead to the expression for $D_t\na P_{J,v}$, and the details are omitted here.

\section{A priori estimates}
\setcounter{equation}{0}
 
We are going to prove Theorem \ref{main theorem}  in this section.  The energy estimate is firstly proved for $J$, and then we go back to the estimate for $v$ and $\G_t$ using Proposition \ref{prop:v}.   
Besides, for the pressure $P$,  since
\[
P=P_{v,v}+\sigma\kappa_\cH,
\]
one can see that the estimate for $P$ relies on the estimates for $v,\G_t$. In a word, the energy estimate for system $\mbox{(WW)}$ means the estimate for $v,\G_t$.

To get started, we recall from the introduction that the energy functional
\[
E(t)=\|\na_{\tau_t}J^\bot\|^2_{L^2(\Gamma_t)} +\|\cD_t J\|^2_{L^2(\Om_t)}+\|\Gamma_t\|^2_{H^{\f52}}+\|v\|^2_{L^2(\Om_t)},
\] 
and the dissipation
\[
F(t)=\big|(\sin \om)\na_{{\bf \tau}_t}J^\perp |_{X_c}\big|^2.
\]

Since the higher-order part in  $E(t)$ and  $F(t)$ relates with $J$,  we need to consider the relationship between $v,\, \G_t$ and $E(t),\, F(t)$. The following proposition shows that $v, \,n_t$ can be controlled by $E(t)$.
\begin{prop}\label{prop:v}
{\it Assuming $\Gamma_t\in H^4$, $v\in H^3(\Om_t)$ and  $\om\in (0,\f \pi 6)$,  one has
\ben\label{est:v}
\|v\|_{H^3(\Om_t)}\leq P(E(t))
\een
and
\ben\label{est:kappa}
\|n_t\|_{H^3(\Gamma_t)}\leq P(E(t)),
\een
where $P(E(t))$ is some polynomial of $E(t)$ with positive constant coefficients.
}
\end{prop}
 
\no {\bf Proof.}
Step 1: estimates for $n_t$. To estimate $n_t$, we will start with $\kappa_\cH$. Since 
\[
\kappa=tr\Pi=\na_{\tau_t}n_t\cdot \tau_t\quad \hbox{with}\quad \na_{\tau_t}n_t\parallel \tau_t\quad\hbox{on}\quad \G_t,
\]
one can see that the higher-order estimate of $n_t$ relies on $\kappa$ or equivalently on $\kappa_\cH$.

In fact, $\kappa_\cH$ satisfies the system below if we have $\cN \kappa=\na_{n_t}\kappa_\cH|_{\G_t}\in H^1(\Gamma_t)$ (which will be proved in the following lines):
\[
\left\{\begin{array}{ll}
\D \kappa_\cH=0,\qquad\hbox{on}\quad\Om_t\\
\na_{n_t}\kappa_\cH|_{\G_t}\in H^1(\Gamma_t),\qquad \na_{n_b} \kappa_\cH|_{\G_b}=0.
\end{array}\right.
\] 
Applying Theorem \ref{est:elliptic-N} and noticing that $\na_{n_t}\kappa_{\cH}|_{\G_t}=J^\perp$, we find
\ben\label{est:kappa1}
\|\kappa_{\cH}\|_{H^{\f52}(\Om_t)}&\leq& C(\|\Gamma_t\|_{H^{\f52}})\|\na_{n_t}\kappa_\cH\|_{H^1(\Gamma_t)}\nonumber\\
&\leq& C(\|\Gamma_t\|_{H^{\f52}})( \|\na_{\tau_t} J^\bot\|_{L^2(\Gamma_t)}+\|\na_{n_t}\kappa_\cH\|_{L^2(\Gamma_t)}).
\een
Using Theorem \ref{trace thm PG} and interpolating as in smooth domains, we have for some constant $\delta>0$, it holds
\ben\label{est:kappa2}
\|\pa_{n_t}\kappa_\cH\|_{L^2(\Gamma_t)}\leq \f{\delta}{ C(\|\Gamma_t\|_{H^{\f52}})} \|\kappa_{\cH}\|_{H^{\f52}(\Om_t)}+\f{ C(\|\Gamma_t\|_{H^{\f52}})}\delta  \|\kappa_{\cH}\|_{H^1(\Om_t)}.
\een
Plugging  \eqref{est:kappa2} into \eqref{est:kappa1} and  choosing a suitable $\delta$ we obtain
\ben\label{est:kappa-1}
\|\kappa_{\cH}\|_{H^{\f52}(\Om_t)}&\leq& C(\|\Gamma_t\|_{H^{\f52}}) ( \|\na_{\tau_t} J^\bot\|_{L^2(\Gamma_t)}+\| \kappa_\cH\|_{H^1(\Om_t)}).
\een
Moreover, we derive from Lemma \ref{Harmonic extension H1 estimate} that
\beno
\| \kappa_\cH\|_{H^1(\Om_t)}\leq \|\Gamma_t\|_{H^{\f52}}
\quad\hbox{and}\quad
\|n_t\|_{H^3(\Gamma_t)}\leq C( \|\kappa\|_{H^2(\Gamma_t)}+\|n_t\|_{L
^2(\Gamma_t)})
\eeno
which together with \eqref{est:kappa-1} imply that
\[
\|\kappa_{\cH}\|_{H^{\f52}(\Om_t)}\le P(E(t)).
\]
As a result, we arrive at by Theorem \ref{trace thm PG}
\beq\label{J estimate}
\|n_t\|_{H^3(\Gamma_t)}\leq P(E(t))\quad\hbox{and}\quad \|J\|_{H^\f32(\Om_t)}\le P(E(t)).
\eeq
Moreover, since 
\[
D_t J=\cD_t J+\na P_{J,v},
\]
we can also derive the estimate for $D_t J$ from \eqref{J estimate}:
\beq\label{D_t J estimate}
\|D_t J\|_{L^2(\Om_t)}\le \|\cD_t J\|_{L^2(\Om_t)}+\|P_{J,v}\|_{H^2(\Om_t)}\le P(E(t))(1+\|v\|_{H^\f52(\Om_t)}),
\eeq
where we applied Theorem \ref{est:elliptic} on $\|P_{J,v}\|_{H^2(\Om_t)}$.

\bigskip

\no Step 2: estimates for $v$.  Recall from $\mbox{(WW)}$ that the velocity $v$ satisfies
\beq\label{condition for v}
\left\{\begin{array}{ll}
\D v=0\qquad\hbox{on}\quad \Om_t,\\
v\cdot n_b|_{\G_b}=0
\end{array}\right.
\eeq
In order to have estimates for $v$, a natural way is to find some condition on $\G_t$ for $v$. We choose  to consider Neumann condition on $\G_t$. As long as we establish the estimate for this condition, we can finish the estimate for $v$ based on the  system above. 

To begin with, we denote the outward normal derivative of $v$ on $\G_t$ by 
\[
\nu=\na_{n_t}v.
\]
For the  estimate  of $\nu$,  we plan to deal with its two components $\nu^\top$ and $\nu^\perp$ respectively.  

\medskip

\no (i) The estimate for $\nu^\top$. This part follows the proof of Proposition 4.3 \cite{SZ} and hence some details are omitted. By a direct computation, we have in our case that
\beno
\Delta_{\Gamma_t}\nu^\top=\na_{\tau_t}(\cD\cdot \nu^\top)\tau_t,
\eeno
which results in the following estimate
 \ben\label{est:nu-1}
 \|\nu^\top\|_{H^{\f32}(\Gamma_t)}\leq C\big(\|\Gamma_t\|_{H^{\f52}}) \|\cD \cdot \nu^{\top}\|_{H^{\f12}(\Gamma_t)}.
 \een
 Next, we give the estimates of $\cD \cdot \nu^{\top}$. From the definition of $\nu$ and a direct computation, we can find   
 \beno
 \cD \cdot \nu^{\top}=\cD_{\tau_t}\nu^\top\cdot\tau_t=\Delta_{\Gamma_t} v\cdot n_t+\Pi(\tau_t)\cdot \na_{\tau_t} v\qquad\hbox{on}\quad \G_t
 \eeno
 where  the second term of the  right side can be controlled by
 \beq\label{est:D_G_t v-R}
 \begin{split}
 \|\Pi(\tau_t)\cdot \na_{\tau_t} v\|_{H^{\f12}(\Gamma_t)}& \le
 C(\|\Gamma_t\|_{H^{\f52}}) \|n_t\|_{H^\f32(\G_t)}\|v\|_{H^{\f52}(\Om_t)}\\
 & \le
C(\|\Gamma_t\|_{H^{\f52}}) \|v\|_{H^{\f52}(\Om_t)}.
\end{split}
 \eeq

On the other hand, for the first term $\Delta_{\Gamma_t} v\cdot n_t$,  recalling equation \eqref{Dt J} for $D_t J$ to  have the estimate
 \beno
 \|D_t J+\na \cH(\Delta_{\Gamma_t} v\cdot n_t)\|_{L^2(\Om_t)}\leq P(E(t))\|v\|_{H^{\f52}(\Om_t)}
 \eeno
where Lemma \ref{Harmonic extension H1 estimate}, Theorem \ref{est:elliptic} (for the $\D^{-1}$ term with $s=2$) and \eqref{J estimate} are applied.
Now moving $D_tJ$ to the right side of the inequality above and applying \eqref{D_t J estimate} and the Poincar\'e inequality lead to 
 \ben\label{est:D_G_t v}
 \|  \cH(\Delta_{\Gamma_t} v\cdot n_t)\|_{H^1(\Om_t)}\leq P(E(t))(1+\|v\|_{H^{\f52}(\Om_t)}).
 \een
 
As a result,  combining  \eqref{est:D_G_t v-R} and \eqref{est:D_G_t v} and  using Theorem \ref{trace thm PG} to get 
 \beno
 \|  \cD \cdot \nu^\top\|_{H^\f12(\Gamma_t)}\leq P(E(t))(1+\|v\|_{H^{\f52}(\Om_t)}),
 \eeno
and furthermore plugging the estimate above into \eqref{est:nu-1}, we finally obtain
 \ben\label{est:nu-top}
 \|\nu^\top\|_{H^{\f32}(\Gamma_t)}\leq P(E(t))(1+\|v\|_{H^{\f52}(\Om_t)}).
 \een

\no (ii) The estimate for $\nu^\perp$.   Since 
\beq\label{relation of nu and v perp}
\nu^\perp=\na_{n_t}v\cdot n_t=\na_{n_t}(v\cdot\tilde n_t)-v\cdot \na_{n_t}\tilde n_t\quad\hbox{on}\quad \G_t
\eeq
for some extension $\tilde n_t$ of $n_t$ on $\Om_t$,  the idea in this part is to deal with $v^\perp:=v\cdot \tilde n_t$ first and then go back to $\nu^\perp$ (the extension $\tilde n_t$ will be defined later). Again, the estimate for $v^\perp$ relies on some elliptic system with boundary conditions.

Firstly, we look for the estimate for $v^\perp$ on $\G_t$. In fact, recalling equation \eqref{Dt k 1} for $\kappa$ to obtain
 \ben\label{est:low-kappa}
 \|D_t \kappa\|_{L^2(\Gamma_t)} \leq C(\|\Gamma_t\|_{H^{\f52}})\|v\|_{H^{\f52}(\Om_t)}
 \een
 and
 \beno
 \|\Delta_{\Gamma_t} v^\bot\|_{H^{\f12}(\Gamma_t)}\leq C(\|\Gamma_t\|_{H^{\f52}}) (\|D_t\kappa\|_{H^{\f12}(\Gamma_t)}+\|v\|_{H^{\f52}(\Om_t)} )
 \eeno
 which implies that
  \ben\label{est:v_bot1}
 \|  v^\bot\|_{H^{\f52}(\Gamma_t)}\leq C(\|\Gamma_t\|_{H^{\f52}}) (\|D_t\kappa\|_{H^{\f12}(\Gamma_t)}+\|v\|_{H^{\f52}(\Om_t)} ).
 \een
 
 To close the estimate above,  we need to find the estimate for $\|D_t\kappa\|_{H^\f12(\G_t)}$.  In fact, we notice that
 \beno
 D_t J=D_t\na \kappa_\cH=\na D_t\kappa_\cH -(\na v)^*\na \kappa_\cH
 \eeno
 which leads to
 \beno
 \|\na D_t\kappa_\cH \|_{L^2(\Om_t)}\leq C(\|(\na v)^* \na \kappa_\cH\|_{L^2(\Om_t)}+\| D_tJ\|_{L^2(\Om_t)}).
 \eeno
On the other hand, for any two points $X_1\in \Om_t$ and $X_2\in \G_t$ where $X_1,\,X_2$ stay in  a vertical line, we can write
\[
D_t\kappa_\cH(X_1)=D_t\kappa(X_2)+\int^{X_1}_{X_2}d D_t\kappa_\cH,
\]
so combining the estimate above for $\|\na D_t\kappa_\cH \|_{L^2(\Om_t)}$ and remembering that our domain $\Om_t$ has a finite depth,  we arrive at 
 \ben\label{est:kappa_trace}
 \|  D_t\kappa_\cH \|_{H^1(\Om_t)}\leq C(\|(\na v)^*\na \kappa_\cH\|_{L^2(\Om_t)}+\| D_tJ\|_{L^2(\Om_t)}+ \|D_t\kappa \|_{L^2(\Gamma_t)}).
 \een
Using Theorem \ref{trace thm PG}, \eqref{D_t J estimate} and the estimate for $\kappa_\cH$ in Step 1 one can  conclude from \eqref{est:kappa_trace} that
   \beno
 \|  D_t\kappa \|_{H^{\f12}(\Gamma_t)}\leq P(E(t))(1 +\|v\|_{H^{\f52}(\Om_t)}+ \|D_t\kappa \|_{L^2(\Gamma_t)}) .
 \eeno
Consequently, plugging the  estimate above and \eqref{est:low-kappa} into \eqref{est:v_bot1} to obtain the estimate for $v^\perp$ on $\G_t$:
  \ben\label{est:v_bot}
 \|  v^\bot\|_{H^{\f52}(\Gamma_t)}\leq P(E(t)) (1+\|v\|_{H^{\f52}(\Om_t)}) .
 \een

Next,  we need also to find the condition for $v^\perp$ on the bottom, $\G_b$. 
To begin with, we define the extension of $n_t$. In fact, one can define Dirichlet-boundary conditions for $\tilde n_t$ as 
\[
\tilde n_t|_{\G_t}=n_t\qquad\hbox{and}\  \tilde n_t=-\tau_b \  \hbox{away from $X_c$ on $\G_b$}.
\] 
Apply a Dirichlet-type trace theorem (for example Theorem 4.7 \cite{MW}) one obtains $\tilde n_t$ defined on $\Om_t$ satisfying
\[
\|\tilde n_t\|_{H^3(\Om_t)}\le C(\|\G_t\|_{H^\f52})\big(\|n_t\|_{H^\f52(\G_t)}+\|\tau_b\|_{H^\f52(\G_b)}\big).
\]

Now we can use  the following condition
\[
\na\cdot v=0\quad\hbox{and}\quad \na\times v=0
\] 
to  find directly that
\[
\begin{split}
\na_{n_b}v^\perp&=\na_{n_b}\Big((\tilde n_t\cdot \tilde\tau_b)v\cdot\tilde\tau_b+(\tilde n_t\cdot\tilde n_b)v\cdot\tilde n_b\Big)\\
&=-(\tilde n_t\cdot n_b)\na_{\tau_b}(v\cdot\tau_b)+r_\perp \qquad\hbox{on}\quad \G_b,
\end{split}
\]
with
\[
\begin{split}
r_\perp=& -(\tilde n_t\cdot \tilde \tau_b)v\cdot(\na_{\tau_b}n_b-\na_{n_b}\tilde\tau_b)+(\tilde n_t\cdot n_b)v\cdot(\na_{\tau_b}\tau_b+\na_{n_b}\tilde n_b)\\
&\quad +\na_{n_b}(\tilde n_t\cdot \tilde \tau_b) v\cdot \tau_b+\na_{n_b}(\tilde n_t\cdot \tilde n_b)v\cdot n_b\qquad\hbox{on}\quad \G_b
\end{split}
\]
where $\tilde n_b,\,\tilde \tau_b$ are  unit orthogonal extensions for $n_b,\,\tau_b$.
On the other hand, since $v\cdot n_b|_{\G_b}=0$, we can arrive at the following oblique condition 
\[
\na_{n_b}v^\perp+b_\perp \na_{\tau_b}v^\perp=R_\perp \qquad\hbox{on}\quad \G_b,
\]
with
\[
b_\perp=\f{\tilde n_t\cdot n_b}{\tilde n_t\cdot \tau_b},\quad R_\perp=r_\perp-b_\perp\na_{\tau_b}(\tilde n_t\cdot \tau_b)\, v\cdot \tau_b.
\]
Summing up the boundary conditions above, we finally conclude the following elliptic system for $v^\perp$:
\[
\left\{\begin{array}{ll}
\D v^\bot =2 \na v\cdot \na \tilde n_t+v\cdot \D \tilde n_t\qquad\hbox{on}\quad\Om_t,\\
 v^\bot  |_{\G_t}\in H^{\f52}(\Gamma_t),
\quad\na_{n_b}v^\bot +b_\perp\na_{\tau_b}v^\bot|_{\G_b} =R_\perp|_{\G_b}.
\end{array}\right.
\]

As a result, based on the  system above, we are ready to find the estimate for $v^\perp$. Since the coefficient $b_\perp$ is a function with bounded support near the contact point, Lemma \ref{est:elliptic-1} is applied here to  have 
\beno
 \|  v^\bot\|_{H^{3}(\Om_t)}\leq  P(E(t))  (1+\|v\|_{H^{\f52}(\Om_t)}),
\eeno
which together with \eqref{relation of nu and v perp}  implies that
\ben\label{est:nu-bot}
 \|  \nu^\bot\|_{H^{\f32}(\Gamma_t)}\leq P(E(t)) (1+\|v\|_{H^{\f52}(\Om_t)}).
\een

\no (iii) The estimate for $v$.
Firstly, combing the estimates \eqref{est:nu-top} and \eqref{est:nu-bot} one can conclude that
\ben\label{est:nu}
 \|\na_{n_t}v\|_{H^\f32(\G_t)}=\|  \nu\|_{H^{\f32}(\Gamma_t)}\leq  P(E(t))(1+\|v\|_{H^{\f52}(\Om_t)}).
 \een

In order to prove the estimate for $v$, we need to check again on system \eqref{condition for v}. We can see that the condition on $\G_b$ is $v\cdot n_b=0$, which is the condition for $v\cdot n_b$ instead of $v$. Therefore, we will derive estimates for $v\cdot\tilde n_b$ and then $v\cdot \tilde \tau_b$ to close the estimate.

In fact, we write down directly the system for $v\cdot \tilde n_b$ as 
\[
\left\{\begin{array}{ll}
\D (v\cdot \tilde n_b)=[\D, \tilde n_b]v\qquad\hbox{on}\quad\Om_t,\\
\na_{n_t}(v\cdot \tilde n_b)|_{\G_t}\in H^{\f32}(\Gamma_t),\qquad v\cdot n_b|_{\G_b}=0.
\end{array}\right.
\] 
where the Neumann condition on the free surface is obtained easily from \eqref{est:nu}. Besides, one can see that this system is again a mixed-boundary problem, although it is slightly different from  system \eqref{eq:elliptic} in Section 5 by switching the two boundary conditions. The elliptic estimate for this system turns out to be  similar as in Theorem \ref{est:elliptic}.  As a result, we arrive at  
\beno
\|v\cdot \tilde n_b\|_{H^{3}(\Om_t)}&\leq &C(\|\Gamma_t\|_{H^{\f52}})(\|[\D, \tilde n_b]v\|_{H^1(\Om_t)}+\|\na_{n_t}(v\cdot \tilde n_b)\|_{H^{\f32}(\Gamma_t)})\\
&\leq & P(E(t)) \big(1+\|v\|_{H^{\f52}(\Om_t)}\big)
\eeno
where we used \eqref{est:nu}.

Next, we will retrieve $v\cdot\tilde \tau_b$ from $v\cdot \tilde n_b$. By the divergence-free and the curl-free conditions for $v$, it is straightforward to show
\[
\begin{split}
\na_{\tilde \tau_b}(v\cdot \tilde \tau_b)&=-\na_{\tilde n_b}(v\cdot \tilde n_b)+v\cdot(\na_{\tilde \tau_b}\tilde \tau_b+\na_{\tilde n_b}\tilde n_b)\qquad \hbox{and}\\
\na_{\tilde n_b}(v\cdot \tilde \tau_b)&=\na_{\tilde\tau_b}(v\cdot \tilde n_b)-v\cdot (\na_{\tilde \tau_b}\tilde n_b-\na_{\tilde n_b}\tilde \tau_b).
\end{split}
\]
Therefore we can have
\[
\|v\cdot \tilde \tau_b\|_{H^{3}(\Om_t)}
\le P(E(t))    (1+\|v\|_{H^{\f52}(\Om_t)}),
\]
which implies that
\[
\|v \|_{H^{3}(\Om_t)}
\le P(E(t))    (1+\|v\|_{H^{\f52}(\Om_t)}).
\]
Finally,  we can conclude the desired estimate by an interpolation.

\ef
 
 \medskip

\subsection{Estimates for $P_{w,v}$ terms}
Before considering about the energy estimates, some more preparations are still needed.  In this part, we present some estimates related to $P_{w,v}$, which will be used frequently later.

In fact, for some vector field $w$ on $\Om_t$, we firstly recall the system for $P_{w,v}$:
 \beno
\left\{
\begin{array}{l}
\Delta P_{w,v}=-  tr(\na w \na v) \qquad  \hbox{on}\quad    \Om_t\\
P_{w,v}|_{\Gamma_t}=0,\quad
\pa_{n_b} P_{w,v} |_{\Gamma_b}=w\cdot \na_v n_b.
\end{array}
\right.
\eeno
The following proposition  gives an elliptic estimate for $P_{w,v}$.
\begin{prop}\label{prop:P(f,v)}
Let  $\om\in(0,\f \pi 6)$ and $s\in[2,4]$. The estimate below holds for any vector field $w\in H^{s-1}(\Om_t)$:
\beno
\|    P_{w, v}\|_{H^{s}(\Om_t)}\leq  P(E(t)) \|w\|_{H^{s-1 }(\Om_t)}.
\eeno
Moreover, if $w=D_t J$, we have 
\beno
\|  P_{D_tJ, v}\|_{H^1(\Om_t)}\leq  P(E(t)) .
\eeno
\end{prop}
\no {\bf Proof.}
The first result comes directly from  Theorem \ref{est:elliptic} (notice that when $s=4$, one needs to apply \eqref{est:kappa} for the coefficient there). So we focus on the second one, which is a variational estimate. 

In fact,  For all $\phi\in \cV=\{\phi\in H^1(\Om_t)\,| \,\phi|_{\Gamma_t}=0\}$, we have the variation equation 
\beq\label{var equation for P DtJ v}
\int_{\Om_t} \na P_{D_tJ, v} \cdot \na \phi dX=\int_{\Om_t} tr(\na D_tJ\na v) \,\phi dX+\int_{\Gamma_b}(D_t J\cdot \na_v n_b )\,\phi ds.
\eeq
  
Firstly,  for the second term on the right side, we  write 
\[
D_t J\cdot \na_v n_b =g\na_{\tau_b} D_t\kappa_\cH-\na_v n_b \cdot (\na v)^* J
\]
where  we notice that $\na_v n_b=(\na_v n_b\cdot \tau_b)\tau_b$ and denote $g=\na_v n_b\cdot \tau_b$. Therefore we obtain 
\beno
\int_{\Gamma_b}(D_t J\cdot \na_v n_b  )\,\phi ds& =& \int_{\Gamma_b}(g\na_{\tau_b} D_t\kappa_\cH-\na_v n_b \cdot (\na v)^* J)\phi ds \\
&\leq& C\|\na_{\tau_b} D_t\kappa_\cH\|_{\wt{H}^{-\f12}(\Gamma_b)}\|g\,\phi\|_{\tilde H^{\f12}(\Gamma_b)}+\|\na_v n_b \cdot (\na v)^* J\|_{L^2(\Gamma_b)}\| \phi\|_{L^2(\Gamma_b)}\\
&\leq&  P(E(t)) \| \phi\|_{H^{1}(\Om_t)},
\eeno
where we used Lemma \ref{H1 zero trace} for $g\,\phi\in\tilde H^\f12(\G_b)$, and moreover, Lemma \ref{G_b H-1/2}, Theorem \ref{trace thm PG}, Proposition \ref{prop:v} and \eqref{est:kappa_trace} for $D_t\kappa_\cH$. Besides, we note that the estimate for $J$ here is given in \eqref{J estimate}.

Secondly, integrating by parts and handling the boundary term similarly as above,  we derive the following estimate for the first term on the right side of \eqref{var equation for P DtJ v}:
\[
\int_{\Om_t} tr(\na D_tJ\na v)\cdot \phi dx 
\le   P(E(t))  \| \phi\|_{H^1(\Om_t)}.
\]

Combining all the estimates above, we conclude that
\[
\int_{\Om_t} \na P_{D_tJ, v} \cdot \na \phi dX\le  P(E(t))   \|w\|_{H^1(\Om_t)},
\]
which leads to the unique existence of $P_{D_tJ,v}\in \cV$ and the desired estimate by Lax-Milgram theorem.

\ef

Based on the estimate for $P_{w,v}$, we will consider  here  some more related terms which will be used later. Moreover,  although the estimates for $P_{J,v}$ and $D_t J$ are already mentioned before,  we still write them down again. 

\medskip

\no 1. $P_{J,v}$ and $P_{v,v}$.  Applying  Proposition \ref{prop:P(f,v)} together with Proposition \ref{prop:v} and \eqref{J estimate},  one has that
\ben\label{est:Pvv}
\|P_{J, v}\|_{H^2(\Om_t)}+\|P_{v,v}\|_{H^{4}(\Om_t)}\leq  P(E(t)) .
\een
\no 2. $H^\f52$ estimate for $P_{J,v}$.   Recalling the definition of $P_{J,v}$ and  applying Theorem \ref{est:elliptic}  to obtain
  \beno
 \| P_{J,v}\|_{H^2(\Om_t)}\leq P(E(t))\|J\|_{H^1(\Om_t)},
 \eeno
 and
 \beno
  \| P_{J,v}\|_{H^3(\Om_t)}\leq  P(E(t)) \|J\|_{H^2(\Om_t)}
 \eeno if we have $J\in H^2(\Om_t)$. 
 Thus, by the complex interpolation theory, we derive
    \beno
  \| P_{J,v}\|_{H^{\f 52}(\Om_t)} \leq   P(E(t))  \| J\|_{H^{\f32}(\Om_t)}.
 \eeno
As a result, combining \eqref{J estimate} we conclude that 
 \ben\label{est:P_{J,v}}
  \| P_{J,v}\|_{H^{\f 52}(\Om_t)}\leq  P(E(t))  .
 \een
\no 3. $D_t J$. Combining \eqref{D_t J estimate} and Proposition \ref{prop:v},   we have
 \beq\label{equ:Relation}
\| D_t J\|_{L^2(\Om_t)}\le   P(E(t)).
 \eeq
\no 4. $D_t P_{v,v}$ and $D_t P_{J,v}$. Recalling that $P_{v,v}=\D^{-1}(-tr(\na v)^2,\,v\cdot \na_vn_b)$, direct computations using \eqref{commutator Dt D-1} lead to 
\[
\begin{split}
D_t P_{v,v}=&\D^{-1}\Big(-tr\big[\big(\na D_tv-(\na v)^*\na v\big)\na v+\na v \big(\na D_t v-(\na v)^*\na v\big)\big],\\
&\quad  \na_vn_b\cdot D_tv+\na_{D_t v}n_b\cdot v-v\cdot(\na v)^*\na n_b\cdot v-\na_v\big((\na v)^*n_b\big)^\top\cdot v\Big)\\
&\quad+\D^{-1}\Big(2\na v\cdot \na^2P_{v,v}+\D v\cdot\na P_{v,v},\,(\na_{n_b}v-\na_vn_b)\cdot \na P_{v,v}\Big)
\end{split}
\] Consequently, applying the Euler equation from $\mbox{(WW)}$, Theorem \ref{est:elliptic}, \eqref{est:Pvv} and a complex interpolation, we can have 
\beq\label{est:D_tPvv}
\|D_t P_{v,v}\|_{H^\f52(\Om_t)}\le P(E(t)).
\eeq 
Similarly, we can also show by a variational argument as in the proof of Proposition \ref{prop:P(f,v)}  that
\beq\label{est:D_t P_{J,v}}
\|D_t P_{J,v}\|_{H^1(\Om_t)}\le P(E(t)).
\eeq

\medskip

\subsection{The dissipation at the contact point}
We now deal with a dissipation term of $J$ on the contact point, which will appear later in the estimates and plays a key role. Besides, one will see from this lemma that a lower bound for the contact angle $\om$ is needed, which can be verified due to the initial contact angle $\omega_0$ and the small time interval $[0,T_0]$ (see the end of this section).
 \begin{lemma}\label{lem:cc}
{\it We have the following equation at the contact point $X_c$
 \[
(D_tJ)^\perp\,( \na_{{\bf \tau}_t}J)^\perp \big|_{X_c}=-\f{\sigma^2}{\beta_c}F(t) +r_c,
\]
where
\beno
r_c= \big(-r +\cot \om(J \cdot  D_t n_b)\big)\, (\sin\om) \na_{{\bf \tau}_t}J^\perp \qquad \textrm{at}\quad X_c
\eeno
with
\[
\begin{split}
r&=-\sigma  \sin\om (\na_{\tau_t}v\cdot D_t n_t)+\sigma  \tau_b\cdot D_t n_t\,(\na_{\tau_t}v\cdot n_t)\\
&\quad+\sigma  \sin\om  \big(\na_{\tau_t}  \na P_{v,v}     -[D_t, \na_{\tau_t} ] v\big)\cdot n_t
- \beta_c D_t \na P_{v,v}\cdot \tau_b.
\end{split}
\]
Moreover if the contact angle $\om\in (0,\f\pi 6)$ and $\sin\om\ge c_0$ for some small constant $c_0>0$, we have the following estimate
\ben\label{est:remaider}
|r_c |\leq  P(E(t))  \,F(t)^{\f12}. 
\een
} 
 \end{lemma}
 \no{\bf Proof.}\,
To begin with,  we recall the boundary condition on $X_c$ from $\mbox{(WW)}$:
\beno
\beta_c v_c=[\gamma]-\sigma \cos \om \qquad \textrm{at}\quad X_c.
\eeno
Moreover, we also have at $X_c$ that
\[
\cos\om=-\tau_t\cdot \tau_b,\quad v_c=-v\cdot\tau_b,
\]
which implies 
\ben\label{equ:cp}
-\beta_c v\cdot \tau_b=[\gamma]+\sigma \tau_t\cdot \tau_b \qquad \textrm{at}\qquad X_c.
\een

On the other hand, we  recall the computations  
\[
D_t n_t=-((\na v)^* n_t)^\top\quad\hbox{and}\quad D_t \tau_t=(\na_{\tau_t}v\cdot n_t)n_t
\qquad \textrm{on}\quad \Gamma_t,
\] 
which will be used in the following lines.

Take $\p_t$ on both sides of \eqref{equ:cp} to obtain
\beno
 -\beta_c (D_t v)\cdot \tau_b=\sigma (D_t \tau_t) \cdot \tau_b\qquad \textrm{at}\quad X_c,
\eeno
where one notices again that $\tau_b, n_b$ are constant vectors near $X_c$. Substituting the Euler equation from $\mbox{(WW)}$ into the  equation above, one gets
\beno
 \beta_c ( \sigma J+\na P_{v,v} + g )\cdot \tau_b=\sigma (D_t \tau_t) \cdot \tau_b\qquad \textrm{at}\quad X_c.
\eeno
Now take $\p_t$ again on  both sides of the equation above to have at $X_c$ that
\beq\label{equ:1}
\begin{split}
\beta_c \sigma D_t J \cdot \tau_b+ \beta_c D_t \na P_{v,v}\cdot \tau_b&=\sigma D_t^2 \tau_t \cdot \tau_b=\sigma  \tau_b\cdot D_t\big( (\na_{\tau_t}v\cdot n_t)\,n_t\big)\\
&=\sigma  \tau_b\cdot (D_t\na_{\tau_t}v\cdot n_t) \,n_t +\sigma  \tau_b\cdot (\na_{\tau_t}v\cdot D_t n_t) \,n_t \\
&\quad+\sigma  \tau_b\cdot (\na_{\tau_t}v\cdot n_t)D_t n_t \\
&=-\sigma  \sin\om   \big(D_t\na_{\tau_t}v\cdot n_t+\na_{\tau_t}v\cdot D_t n_t\big )\\
&\quad+\sigma  \tau_b\cdot D_t n_t\,(\na_{\tau_t}v\cdot n_t),
\end{split}
\eeq
where we used 
\[ 
n_t\cdot \tau_b=-\sin\om\qquad\hbox{at}\quad X_c.
\]

To deal with the terms on the right side, apply $\na_{\tau_t}$ on the Euler equation in $\mbox{(WW)}$ (constrained on $\G_t$) to arrive at
\beno
D_t \na_{\tau_t }v\cdot n_t= (\na_{\tau_t} D_t v)\cdot n_t+[D_t, \na_{\tau_t} ] v\cdot n_t=-\na_{\tau_t}(\sigma J  + \na P_{v,v} )\cdot n_t+[D_t, \na_{\tau_t} ] v\cdot n_t,
\eeno
which can be substituted into \eqref{equ:1} to obtain
\beq\label{equ:J1}
\beta_c \sigma D_t J \cdot \tau_b= \sigma^2(\sin \om)\na_{\tau_t}J \cdot n_t+r \qquad \textrm{at}\quad X_c,  
\eeq
where $r$ is defined by
\[
\begin{split}
r=&\sigma  \sin\om  \big(\na_{\tau_t}  \na P_{v,v}     -[D_t, \na_{\tau_t} ] v\big)\cdot n_t
-\sigma  \sin\om (\na_{\tau_t}v\cdot D_t n_t)\\
&\quad+\sigma  \tau_b\cdot D_t n_t\,(\na_{\tau_t}v\cdot n_t)-\beta_c D_t \na P_{v,v}\cdot \tau_b.
\end{split}
\]

To retrieve $D_t J$, we still need to compute  $D_t J\cdot n_b$. In fact, we recall the definition of $J$ to find 
\beno
J\cdot n_b=0\qquad \textrm{on}\quad \Gamma_b.
\eeno
Taking $\p_t$ on  both sides of the condition above leads to
\ben\label{equ:J2}
D_t J \cdot n_b=- J \cdot D_t n_b\qquad \textrm{on}\quad \Gamma_b.
\een

Combining \eqref{equ:J1} and  \eqref{equ:J2} to obtain that
\beq\label{equ:cc}
 D_tJ = \f{\sigma^2 }{\beta_c}\sin\om (\na_{\tau_t} J)^\bot 
 \tau_b+\big(r \,\tau_b-  ( J \cdot  D_t n_b)\,n_b\big) \qquad \textrm{on}\qquad \Gamma_b.
\eeq
As a result,  we have 
 \[
(D_tJ)^\perp\, (\na_{{\bf \tau}_t}J)^\perp|_{X_c}=-\f{\sigma^2}{\beta_c}  F(t)+r_c.
\]
with 
\[
r_c=\big(r \,\tau_b-  ( J \cdot  D_t n_b)\,n_b\big)\cdot n_t\na_{{\bf \tau}_t}J^\perp.
\]

In the end, we prove the estimate \eqref{est:remaider}. In fact, applying Lemma \ref{embedding}, Proposition \ref{prop:v}, Proposition \ref{prop:P(f,v)} and \eqref{est:D_tPvv}, it's straightforward to prove 
\beno
\big|r|_{X_c}\big|+\big|J \cdot  D_t n_b|_{X_c}\big| \leq P(E(t))
\eeno
with $D_t n_b=-\big((\na v)^*n_b\big)^\top$.
Therefore,  the proof is  finished.

\ef

\bigskip

\subsection{Energy estimates.}
Finally,  it's the time to prove the a priori estimate in forms of $E(t),F(t)$. We will start with \eqref{eqn cDt2 J} here. For the sake of simplicity, we perform the estimates directly on $J$. Although, to be more strict, the following estimates  should be performed firstly on a sequence of  smooth functions converging to $J$, and then we should  show that the final energy estimate also holds when the limit is taken. 

To begin with,  taking the inner product with $\cD_t J$ on both sides of \eqref{eqn cDt2 J} and integrating on $\Om_t$, one  obtains that
\ben\label{equ:Ee2}
 \int_{\Om_t} D_t \cD_t J\cdot\cD_t JdX +\sigma\int_{\Om_t} \cA J\cdot \cD_t JdX=\int_{\Om_t} R\cdot \cD_t JdX
\een
where 
\[
R=R_0+D_t\na P_{J,v}.
\]
For the first term on the left side of \eqref{equ:Ee2}, we simply have
\beno
\int_{\Om_t}D_t  \cD_t J\cdot\cD_t JdX
&=&\f12 \pa_t \int_{\Om_t} |\cD_t J|^2dX .
\eeno
 
Applying \eqref{est:D_t P_{J,v}} on $D_t\na P_{J, v}$  term from the right side, we find 
\beno
\int_{\Om_t}D_t\na P_{J, v} \cdot  \cD_t J dX\leq \|D_t\na P_{J, v}  \|_{L^2(\Om_t)}\|\cD_t J  \|_{L^2(\Om_t)}\leq P(E(t)).
\eeno

Next,  For the second term on the left side of  \eqref{equ:Ee2},  one derives that
\[
\begin{split}
\int_{\Om_t} \cA J\cdot \cD_t JdX&=\int_{\Om_t} \na\mathcal{H}(-\Delta_{\Gamma_t} J^\bot)\cdot \cD_t JdX\\
&=-\int_{\Gamma_t} \Delta_{\Gamma_t} J^\bot \, (\cD_t J\cdot n_t)ds.
\end{split}
\]
Recalling 
\[
\cD_tJ=D_t J+\na P_{J,v}\quad\hbox{and}\quad \cD_t J\circ u\,\in T_{u(t)}\G
\]
 from the Hodge decomposition in Section 3,   one can have by integrating by parts that
\beq\label{equ:Ee1}
\begin{split}
\int_{\Om_t} \cA J\cdot \cD_t J dX
&=-\int_{\Gamma_t} \Delta_{\Gamma_t} J^\bot   (D_t J\cdot n_t) ds-\int_{\Gamma_t} \Delta_{\Gamma_t} J^\bot  (\na P_{J,v}\cdot n_t)ds\\
&=\int_{\Gamma_t} \na_{\tau_t} J^\bot \, \na_{\tau_t}(D_t J\cdot n_t )ds-(D_t J)^\bot \na_{\tau_t} J^\bot \big|_{X_c} -\int_{\Gamma_t} \Delta_{\Gamma_t} J^\bot  (\na  P_{J,v}\cdot n_t)ds.
\end{split}
\eeq
We will check these terms one by one. In fact, for the first term on the right side of \eqref{equ:Ee1}, one deduces that
\[
\begin{split}
&\int_{\Gamma_t} \na_{\tau_t} J^\bot\, \na_{\tau_t}( D_t J\cdot n_t )ds\\
&=\int_{\Gamma_t} \na_{\tau_t} J^\bot    \,D_t( \na_{\tau_t} J^\bot )ds-\int_{\Gamma_t} \na_{\tau_t} J^\bot \,[  D_t, \na_{\tau_t}] J^\bot ds - \int_{\Gamma_t} \na_{\tau_t} J^\bot  \,\na_{\tau_t}(J\cdot  D_t  n_t )ds\\
&=\f12\,\p_t\int_{\Gamma_t}| \na_{\tau_t} J^\bot|^2ds-\int_{\Gamma_t} \na_{\tau_t} J^\bot \, [  D_t, \na_{\tau_t}] J^\bot ds- \int_{\Gamma_t} \na_{\tau_t} J^\bot \, \na_{\tau_t}(J\cdot  D_t  n_t )ds.
\end{split}
\]
For the second term on the right side of \eqref{equ:Ee1},  applying Lemma  \ref{lem:cc}  leads to
\[
-  ( D_t J)^\bot \na_{\tau_t} J^\bot\big|_{X_c}
=\f{\sigma^2}{\beta_c}  F(t)-r_c
\]
where
\[
|r_c |\leq  P(E(t)) \, F(t)^{\f12}
\]
as long as 
\[
\sin\om\ge c_0>0\quad\hbox{for some constant}\, c_0,
\]
which will be checked in the end of this paper.
 
Combining all the  estimates above, we have
\[
\begin{split}
&\partial_t\Big(\f12 \int_{\Om_t}|\cD_t J|^2dX+\f\sigma2 \int_{\Gamma_t}| \na_{\tau_t} J^\bot|^2ds\Big)+ \f{\sigma^3}{2\beta_c}  F(t)\\
&\le P(E(t))+ \sigma\int_{\Gamma_t} \Delta_{\Gamma_t} J^\bot  \na P_{J,v}\cdot n_tds +\sigma \int_{\Gamma_t} \na_{\tau_t} J^\bot \, [ D_t, \na_{\tau_t}] J^\bot ds \\
&\qquad+\sigma\int_{\Gamma_t} \na_{\tau_t} J^\bot \,\na_{\tau_t}(J\cdot D_t  n_t )ds +\int_{\Om_t} R_0\cdot \cD_t JdX.
\end{split}
\]
So now it remains to deal with the right side of the energy estimate above. In fact, it is straightforward to show from Theorem \ref{trace thm PG}, Proposition \ref{prop:v} and  \eqref{J estimate} that 
\[
\int_{\Gamma_t} \na_{\tau_t} J^\bot \, [  D_t, \na_{\tau_t}] J^\bot ds+\int_{\Gamma_t} \na_{\tau_t} J^\bot \,\na_{\tau_t}(J\cdot  D_t  n_t )ds  \leq  P(E(t)). 
\]
Moreover, a direct computation and applying Theorem \ref{trace thm PG}, \eqref{J estimate} and \eqref{est:P_{J,v}} imply 
\[
\begin{split}
 \int_{\Gamma_t} \Delta_{\Gamma_t} J^\bot  \,\na P_{J,v}\cdot n_tds&=-\int_{\Gamma_t} \na_{\tau_t} J^\bot \,\na_{\tau_t} \na P_{J,v}\cdot n_tds+ \big( \na_{\tau_t} J^\bot \, \na P_{J,v}\cdot n_t\big)\big|_{X_c}\\
 &\le \|\na_{\tau_t} J^\bot \|_{L^2(\G_t)} \,\|\na_{\tau_t} \na P_{J,v}\cdot n_t\|_{L^2(\G_t)}+  F(t)^\f12\,\Big|\f1{\sin\om}\na P_{J,v}\cdot n_t|_{X_c}\Big|\\
  &\le \f14 \f{\sigma^2}{\beta_c}  F(t) + P(E(t))
\end{split}
\]
as long as we have $\sin\om\ge c_0>0$.

As a result, we arrive at 
\beno
 \pa_t \Big(\f12\int_{\Om_t}|\cD_t J|^2dX+\f\sigma2\int_{\Gamma_t}| \na_{\tau_t} J^\bot|^2ds\Big)+ \f{\sigma^3}{4\beta_c} F(t) \leq   P(E(t)) + \| R_0\|^2_{L^2(\Om_t)},
\eeno
which tells us  that, to close the energy estimates, the only thing left is to prove the estimate for the reminder term $R_0$. 
  
\medskip
To begin with, one can see already that the boundary conditions on $\G_b$ play an important role in the variational estimates, which are handled differently compared to the smooth-domain case.  Therefore, a lemma is presented here focusing on a typical type of boundary conditions needed in the estimate for $R_0$.
\begin{lemma}\label{na_n_b v-na_v n_b boundary}
{\it 
Let $w\in H^1(\Om_t)$ satisfying $\D w\in L^2(\Om_t)$. Then  the boundary condition
\[
(\na_{n_b}v-\na_v n_b)\cdot \na w\big|_{\G_b}
\] 
 makes sense in the variation formulation: For any $\phi\in \cV=\{\phi\in H^1(\Om_t)\,\big|\,\phi|_{\G_t}=0\}$, one has 
\[
\displaystyle\int_{\G_b}(\na_{n_b}v-\na_v n_b)\cdot \na w\,\,\phi \,ds\le P(E(t))\big(\|w\|_{H^1(\Om_t)}+\|\D w\|_{L^2(\Om_t)}\big)\|\phi\|_{H^1(\Om_t)}.
\]
}
\end{lemma}

\no {\bf Proof.}  The proof  lies in clarifying the space for the boundary condition. 

Firstly, one needs to rewrite $\na_{n_b}v-\na_v n_b$.  Since $n_b$ is the unit normal vector on $\G_b$, one knows that $\na_vn_b=(\na_v n_b\cdot\tau_b)\tau_b$. Moreover, decompose $\na_{n_b}v$  with respect to $\tau_b$ and $n_b$ to have 
\[
\na_{n_b}v-\na_v n_b=(\na_{n_b}v\cdot\tau_b-\na_v n_b\cdot\tau_b)\tau_b+(\na_{n_b}v\cdot n_b)n_b.
\]
Consequently, the boundary condition can be written as 
\[
(\na_{n_b}v-\na_v n_b)\cdot \na w|_{\G_b}=(\na_{n_b}v\cdot\tau_b-\na_v n_b\cdot\tau_b)\na_{\tau_b}w+(\na_{n_b}v\cdot n_b)\na_{n_b} w\big|_{\G_b}.
\]

Secondly, since $w\in H^1(\Om_t)$, one has $w\in H^\f12(\G_b)$ by Theorem \ref{trace thm PG}. Applying Lemma \ref{G_b H-1/2}, one can see also that $\na_{\tau_b}w\in \tilde H^{-\f12}(\G_b)$ with the estimate
\[
\|\na_{\tau_b}w\|_{\tilde H^{-\f12}(\G_b)}\le P(E(t))\|w\|_{H^\f12(\G_b)}\le P(E(t))\|w\|_{H^1(\Om_t)}.
\]
On the other hand, use Lemma \ref{na_n_b w on G_b} to find $\na_{n_b}w|_{\G_b}\in\tilde H^{-\f12}(\G_b)$ satisfying
\[
\|\na_{n_b}w\|_{\tilde H^{-\f12}(\G_b)}\le P(E(t))\big(\|w\|_{H^1(\Om_t)}+\|\D w\|_{L^2(\Om_t)}\big).
\]
Finally, summing up the estimates above leads  to 
\[
\begin{split}
&\displaystyle\int_{\G_b}(\na_{n_b}v-\na_v n_b)\cdot \na w\,\,\phi \,ds\\
 &\quad \le \|\na_{\tau_b}w\|_{\tilde H^{-\f12}(\G_b)}\|(\na_{n_b}v\cdot\tau_b-\na_v n_b\cdot\tau_b)\phi\|_{\tilde H^{\f12}(\G_b)}
+\|\na_{n_b}w\|_{\tilde H^{-\f12}(\G_b)}   \|(\na_{n_b}v\cdot n_b)\phi\|_{\tilde H^{\f12}(\G_b)}\\
&\quad \le P(E(t))\big(\|w\|_{H^1(\Om_t)}+\|\D w\|_{L^2(\Om)}\big)\|\phi\|_{H^1(\Om_t)},
\end{split}
\]
where Theorem \ref{trace thm PG}, Lemma \ref{H1 zero trace}, Lemma \ref{embedding} and Proposition \ref{prop:v} are applied.

\ef

Now it's the time to prove the estimate for $R_0$.
\begin{proposition}
{\it
We have the following estimate for the remainder term $R_0$ defined in \eqref{eqn for Dt2 J}:
\beno
\|R_0\|_{L^2(\Om_t)}\leq  P(E(t)) .
\eeno
}
\end{proposition}
\no {\bf Proof.}
Recall from \eqref{eqn for Dt2 J} that 
\[
R_0=-\s\na \cH(J\cdot \D_{\G_t}n_t)+\na \cH(n_t\cdot \D_{\G_t}\na P_{v,v})+\na\cH(R_1) +A_1+A_2+A_3
\] 
where $R_1$ and $A_1,A_2,A_3$ are defined in \eqref{Dt2 kappa} and \eqref{A_1}, \eqref{A_2}, \eqref{Dt2 J eqn1} respectively, so the estimate for $R_0$ lies in the estimates for all the terms above.

\medskip
\no - Estimate for $\s\na \cH(J\cdot \D_{\G_t}n_t)$.  Applying Theorem \ref{trace thm PG}, Lemma \ref{Harmonic extension H1 estimate} and \eqref{J estimate} one finds that
\[
\|\s\na \cH(J\cdot \D_{\G_t}n_t)\|_{L^2(\Om_t)}\le P(E(t))\|J\cdot \D_{\G_t}n_t\|_{H^{\f12}(\G_t)}\le  P(E(t)).
\]
 
\no - Estimate for $\na \cH(n_t\cdot \D_{\G_t}\na P_{v,v})$. Similarly as the previous term, one has 
\[
\|\na \cH(n_t\cdot \D_{\G_t}\na P_{v,v})\|_{L^2(\Om_t)}\le  P(E(t)).
\]

\no- Estimate for $\na\cH(R_1) $.  Examining the expression for $R_1$ from \eqref{Dt2 kappa}  carefully, one can see that the leading-order terms in $R_1$ are like $\na^2 v$,  $\na n_t$, $\kappa$ and $\na^2P_{v,v}$.  Therefore, $R_1$ can be dealt directly and the details are omitted here.  As a result, we have
\[
\|\na\cH(R_1) \|_{L^2(\Om_t)}\le  P(E(t)).
\]

\no- Estimate for $A_1$. Recall from \eqref{A_1} that 
\[
A_1=\na w
\] with the notation
\[
w= \D^{-1}\big(2\na v\cdot\na^2\cH(D_t\kappa)+\D v\cdot\na \cH(D_t\kappa) ,(\na_{n_b}v-\na_vn_b)\cdot\na \cH(D_t\kappa)\big),
\]
so one needs to deal with the estimate for $w$. In fact, by the notation $\D^{-1}$ we know that $w$ satisfies the system
\[
\left\{\begin{array}{ll}
\D w=2\na v\cdot\na^2\cH(D_t\kappa)+\D v\cdot\na \cH(D_t\kappa),\qquad\hbox{on}\quad \Om_t,\\
w|_{\G_t}=0,\quad \na_{n_b}w|_{\G_b}=(\na_{n_b}v-\na_vn_b)\cdot\na \cH(D_t\kappa)|_{\G_b}
\end{array}
\right.
\]
which is defined by the variation equation 
\[
\begin{split}
\displaystyle\int_{\Om_t}\na w\cdot\na\phi dX=&
-\displaystyle\int_{\Om_t}\Big(2\na v\cdot\na^2\cH(D_t\kappa)+\D v\cdot\na \cH(D_t\kappa)\Big)\, \phi dX\\
&\quad
+\displaystyle\int_{\G_b}(\na_{n_b}v-\na_vn_b)\cdot\na \cH(D_t\kappa)\,\phi ds
\end{split}
\]
for any $\phi \in \cV=\{\phi\in H^1(\Om_t)\,\big|\,\phi|_{\G_t}=0\}$. 

The estimate for $w$ can be derived from the variation equation. To begin with,  one has from the proof (ii) of Proposition \ref{prop:v} that
\[
\|D_t\kappa\|_{H^\f12(\G_t)}\le  P(E(t)),
\]
which together with Lemma \ref{Harmonic extension H1 estimate} leads to
\[
\|\cH(D_t\kappa)\|_{H^1(\Om_t)}\le C\,\|D_t\kappa\|_{H^\f12(\G_t)}\le  P(E(t)).
\]

Secondly, we consider the boundary term in the variation equation. Indeed, one can check directly that $\cH(D_t\kappa)\in H^1(\Om_t)$ satisfies the condition in Lemma \ref{na_n_b v-na_v n_b boundary}, so applying this lemma leads to 
\[
\Big|\displaystyle\int_{\G_b}(\na_{n_b}v-\na_v n_b)\cdot \na \cH(D_t\kappa)\,\,\phi \,ds\Big|
\le   P(E(t)) \|\phi\| _{H^1(\Om_t)}.
\]

Moreover,  integrating by parts, applying Lemma \ref{embedding} and Proposition \ref{prop:v}, and dealing with the boundary term similarly as above, one obtains the estimate for the first term on the right side of the variation equation:
\[
\big|\displaystyle\int_{\Om_t}\Big(2\na v\cdot\na^2\cH(D_t\kappa)+\D v\cdot\na \cH(D_t\kappa)\Big)\, \phi dX\big|
\le  P(E(t)) \|\phi\|_{H^1(\Om_t)}.
\]

Summing up these estimates above, we  conclude by Lax-Milgram Theorem that  the variation equation for $w$ admits a unique solution $w\in H^1(\Om_t)$ with the estimate
\[
\|w\|_{H^1(\Om_t)}\le  P(E(t)) ,
\]
which implies
\[
\|A_1\|_{L^2(\Om_t)}\le  P(E(t)) .
\]

\no- Estimate for $A_2$. The estimate is similar as that for $A_1$. Recall from \eqref{A_2} that 
\[
\begin{split}
A_2&=\na\D^{-1}\big(2\na v\cdot\na^2w_{A2}+\D v\cdot \na w_{A2},(\na_{n_b}v-\na_vn_b)\cdot\na w_{A2}\big) +\na\D^{-1}(h_{A2},\, g_{A2})
\\
&:=A_{21}+A_{22}
\end{split}
\]
where
\[
\begin{split}
 w_{A2}=&\D^{-1}\big(2\na v\cdot \na J+\D v\cdot J,\,(\na_{n_b}v-\na_v n_b)\cdot J\big),\\
 h_{A2}=&2\na v\cdot(\na D_tJ-(\na v)^*J)+2(\na D_tv-(\na v)^*\na v)\cdot \na J+D_tJ\cdot\D v\\
 &\quad+J\cdot (\D D_tv-\D v\cdot\na v-2\na v\cdot \na^2 v),\\
g_{A2}=&(\na_vn_b-\na_{n_b}v)\cdot \na v\cdot J+\na_{n_b}D_t v\cdot J-(D_t v-\na_v v)\cdot \na n_b\cdot J\\
&\quad +\na_v\big((\na v)^*n_b\big)^\top\cdot J+(\na_{n_b}v-\na_vn_b)\cdot D_tJ.
\end{split}
\]

Firstly, in order to deal with $A_{21}$, one needs to handle $w_{A2}$. In fact,  $w_{A2}$ satisfies the system
\[
\left\{\begin{array}{ll}
\D w_{A2}=2\na v\cdot \na J+\D v\cdot J,\qquad \hbox{on}\quad \Om_t\\
w_{A2}|_{\G_t}=0,\quad \na_{n_b}w_{A2}|_{\G_b}=(\na_{n_b}v-\na_v n_b)\cdot J.
\end{array}\right.
\] 
Noticing that $J=\na\kappa_{\cH}\in H^\f32(\Om_t)$, the variational estimate  for $w_{A2}\in H^1(\Om_t)$ can be done similarly as the $w$ system in the estimate for $A_1$, which turns out to be 
\beno
\|w_{A2}\|_{H^1(\Om_t)}
&\le& C\,\big(\|2\na v\cdot \na J+\D v\cdot J\|_{L^2(\Om_t)}
+P(E(t))\|\kappa_{\cH}\|_{H^1(\Om_t)}\big)\\
&\le&   P(E(t)) .
\eeno
Moreover, since 
\[
\D w_{A2}=2\na v\cdot \na J+\D v\cdot J\in L^2(\Om_t),
\]
 we also have that $w_{A2}\in E(\D;L^2(\Om_t))$. 

Now we can close the estimate for $A_{21}$.  In fact, the variational estimate for $A_{21}$ is again similar as before, so we omit the details to write directly that 
\[
\|A_{21}\|_{L^2(\Om_t)}\le P(E(t)),
\]
where Lemma \ref{na_n_b v-na_v n_b boundary} is applied to the boundary condition on $\G_b$.

Secondly, we consider the estimate for $A_{22}$.  Denoting $u_{A2}=\D^{-1}(h_{A2},g_{A2})$, we have 
\[
A_{22}=\na u_{A2}.
\]
Similarly as before, we deal with the variational estimate for $u_{A2}$. In fact, $u_{A2}$ is defined by the variation equation 
\beq\label{var eqn for u_A2}
\displaystyle\int_{\Om_t}\na u_{A2}\cdot\na \phi dX=-\displaystyle\int_{\Om_t} h_{A2}\,\phi dX+\displaystyle\int_{\G_b}g_{A2}\,\phi ds
\eeq
where $\phi\in \cV=\{\phi\in H^1(\Om_t)\,\big|\,\phi|_{\G_t}=0\}$. So the variational estimate again lies in the estimates for the two integrals on the right side.

For the term of $h_{A2}$,  direct estimates as before lead to  
\beq\label{estimate h_A2}
\Big|\displaystyle\int_{\Om_t} h_{A2}\,\phi dX\Big|\le P(E(t))\|\phi\|_{H^1(\Om_t)}.
\eeq

On the other hand, we consider the estimate for $\displaystyle\int_{\G_b}g_{A2}\,\phi ds$ from \eqref{var eqn for u_A2}.  Plugging in the expression for $g_{A2}$ one arrives at
\[
\begin{split}
&\displaystyle\int_{\G_b}g_{A2}\,\phi ds\\
&=
\displaystyle\int_{\G_b}(\na_vn_b-\na_{n_b}v)\cdot \na v\cdot J\,\phi ds+\displaystyle\int_{\G_b}\na_{n_b}D_t v\cdot J\,\phi ds-\displaystyle\int_{\G_b}(D_t v-\na_v v)\cdot \na n_b\cdot J\,\phi ds\\
&\qquad +\displaystyle\int_{\G_b}\na_v\big((\na v)^*n_b\big)^\top\cdot J\,\phi ds+\displaystyle\int_{\G_b}(\na_{n_b}v-\na_vn_b)\cdot D_tJ\,\phi ds\\
&:= B_1+B_2+\dots+B_5,
\end{split}
\]
and the terms will be checked one by one. 

In fact, similar estimates as before imply 
\[
|B_1+\dots+B_4|\le P(E(t))\|w\|_{H^1(\Om_t)}.
\]
It remains to deal with the last term $B_5$. Since
\[
D_tJ=D_t\na\kappa_{\cH}=\na D_t\kappa_{\cH}-(\na v)^*J
\]
where  $D_t\kappa_{\cH}\in H^1(\Om_t)$ and 
\[
\D D_t\kappa_{\cH}=2\na v\cdot J+\D v\cdot J\in L^2(\Om_t). 
\]
As a  result, one has 
\[
\begin{split}
|B_5|&\le \Big|\displaystyle\int_{\G_b}(\na_{n_b}v-\na_v n_b)\cdot \na D_t\kappa_{\cH} \,\phi ds\Big|
+
\Big|\displaystyle\int_{\G_b}(\na_{n_b}v-\na_v n_b)\cdot (\na v)^*J\,\phi ds\Big|\\
&\le P(E(t))\big(\|D_t\kappa_{\cH}\|_{H^1(\Om_t)}+\|\D D_t\kappa_{\cH}\|_{L^2(\Om_t)}\big)\|\phi\|_{H^1(\Om_t)}\\
&\qquad+\|(\na_{n_b}v-\na_v n_b)\cdot (\na v)^*J\|_{L^2(\G_b)}  \|\phi\|_{L^2(\G_b)}\\
&\le P(E(t))\|\phi\|_{H^1(\Om_t)}
\end{split}
\] 
where we applied Lemma \ref{na_n_b v-na_v n_b boundary}, Lemma \ref{embedding}, Proposition \ref{prop:v} and \eqref{est:kappa_trace}.

Summing up the estimates from $B_1$ to $B_5$, we have
\beq\label{estimate g_A2}
\Big|\displaystyle\int_{\G_b}g_{A2}\,\phi ds\Big|\le P(E(t))\|\phi\|_{H^1(\Om_t)}.
\eeq

Combining  \eqref{estimate h_A2} and \eqref{estimate g_A2}  above, we finally conclude that the variation equation \eqref{var eqn for u_A2} admits an unique solution   $u_{A2}\in \cV$ with the estimate 
\[
\|A_{22}\|_{L^2(\Om_t)}=\|\na u_{A2}\|_{L^2(\Om_t)}\le P(E(t)).
\]
As a result, we  finish the estimate for $A_2$. 

\no- Estimate for $A_3$. This term can be handled directly to arrive at
\[
\|A_3\|_{L^2(\Om_t)}\le P(E(t)).
\]
Therefore, summing up all the estimates above, the proof is finished.

\ef

Now, we can finally conclude that 
\beq\label{est:E1}
 \pa_t \Big(\f12\int_{\Om_t}|\cD_t J|^2dX+\f\sigma2\int_{\Gamma_t}| \na_{\tau_t} J^\bot|^2ds\Big)+ \f{\sigma^3}{4\beta_c} F(t)  \leq   P(E(t)) .
\eeq

To close the energy, we still need to deal with $\|\Gamma_t\|_{H^{\f 52}}$ and $\|v\|_{L^2(\Om_t)}$.  Firstly, applying the Euler equation from $\mbox{(WW)}$ with \eqref{est:Pvv} and \eqref{J estimate}, we can prove directly that
\beq\label{est:v_low}
\pa_t\|v\|^2_{L^2(\Om_t)}\leq   P(E(t)) .
\eeq

Secondly, for the term $\|\Gamma_t\|_{H^{\f 52}}$,  we can parametrize $\G_t$ under Eulerian coordinates $(x,z)$ by  
\[
\G_t=\big\{(x,z)\,|\,z=\eta(t,x),\,t>0,\, x\ge c(t)\big\}
\]
where  $c(t)$ is the $x$ coordinate for the contact point $X_c$. Therefore one can write
\beno
\|\Gamma_t\|_{H^{\f52}}= \|\eta\|_{H^{\f52}(c(t),\infty)}
\eeno
and the estimate for $\G_t$ means the estimate for $\eta$.

On the other hand, notice that \mbox{(WW)} contains the kinematic condition on $\G_t$, which can be written in form of $\eta,\,v$ as 
\beno
\pa_t\eta +v_1\, \pa_x\eta=v_2\qquad \textrm{on}\quad [c(t),\infty)
\eeno
where  $v=(v_1,v_2)^t$.  Moreover, a direct computation shows that under this parametrization, the material derivative $D_t$ on $\G_t$  is simply 
\[
D_t=\partial_t+v_1\partial_x,
\]
so the equation above for $\eta$ can be rewritten as 
\[
D_t\eta=v_2.
\]
As a result, it is straightforward  to show that 
 \beq\label{est:eat_low}
 \partial_t \|\G_t\|^2_{H^\f52}=\partial_t\|\eta\|^2_{H^{\f52}(c(t),\infty)}\le P(E(t)).
 \eeq
Summing up  \eqref{est:E1}, \eqref{est:v_low} and \eqref{est:eat_low}, we finally arrive at the energy estimate
\[
\partial_t E(t)+F(t)\le P(E(t)),
\]
and integrating on both sides on a time interval $[0,T_0]$ (to be fixed later) leads to
\[
 \sup_{t\in [0,T_0]} E(t) +  \int_0^{T_0} F(t)    \leq   E(0)+ \int_0^{T_0} P(E(t)) .
\]

In the end, it remains to consider about the evolution of the contact angle $\om(t)$ and verify the condition 
\[
\sin \om(t)\ge c_0>0\quad\hbox{on}\ [0,T_0],\ \hbox{for some constant } c_0.
\] 
In fact, we have at the initial time $t=0$ that 
\[
\om(0)\in (0,\f\pi 6)\quad\hbox{so}\quad \sin\om(0)=-n_t(0)\cdot \tau_b(0)|_{X_c(0)} >0, 
\]
where $\tau_b(0)$ is a constant vector since we set in the beginning that  $\G_b$ becomes a straight line segment near the contact point. So we choose $T_0$ small enough such that for any $t\in [0, T_0]$,  $\tau_b(t)|_{X_c}=\tau_b(0)|_{X_c(0)}$ to  show that
\[
\begin{split}
|\sin \om(0)-\sin\om(t)|&=\big|n_t(t)|_{X_c(t)}-n_t(0)|_{X_c(0)}\big|\,|\tau_b(0)|\\
&\le T_0\sup_{ [0,T_0]}\big|D_t n_t|_{X_c}\big|\,|\tau_b(0)|\\
&\le T_0 \sup_{ [0,T_0]}P(E(t))
\end{split}
\]
Consequently, when $T_0$ is small enough, one finds a small constant $\delta>0$ such that
\[
|\sin \om(0)-\sin\om(t)|<\delta\qquad\hbox{for any }\  t\in [0,T_0],
\]
which infers
\[
\sin\om(t)\ge c_0,\qquad  t\in [0,T_0]
\]
for some constant $c_0>0$. Meanwhile, we can also have 
\[
\om(t)\in (0,\f\pi 6),\qquad t\in [0,T_0].
\]
As a result, our main theorem is proved.

\bigskip

\noindent{\bf Acknowledgement}.  The authors would like to thank Chongchun Zeng for fruitful discussions. The author Mei Ming is supported by NSFC no.11401598.

\end{document}